\documentclass[aap,MSNbibl,nameyear,seceqn,dvips]{arximspdf}
\usepackage{dcolumn}

\doi{10.1214/12-AAP892} 
\volume{23}
\issue{5}
\pubyear{2013}
\firstpage{1988}
\lastpage{2019}

\makeatletter
\newcolumntype{d}[1]{D{.}{.}{#1}}
\newcommand{\rrVert}{\Vert}
\newcommand{\rrvert}{\vert}
\newcommand{\llVert}{\Vert}
\newcommand{\llvert}{\vert}
\newcommand{\eqref}[1]{(\ref{#1})}
\newtheorem{thmm}{Theorem}[section]
\newtheorem{lem}[thmm]{Lemma}
\newtheorem{lemm}{Lemma}[section]
\newtheorem{propp}[lemm]{Proposition}
\newtheorem{cor}[thmm]{Corollary}
\newtheorem{prop}[thmm]{Proposition}
\newproclaim{defi}[thmm]{Definition}
\newproclaim{rem}[thmm]{Remark}
\newcommand{\E}{\mathrm{E}}
\renewcommand{\P}{\mathrm{P}}
\newcommand{\Var}{\operatorname{Var}}

\newcommand{\arginf}{\operatorname{arg\,inf}}
\newcommand{\esssup}{\operatorname{ess\,sup}}
\makeatother

\begin{document}
\begin{frontmatter}

\title{Solving optimal stopping problems via empirical dual optimization}
\runtitle{Empirical dual optimization}

\begin{aug}
\author[A]{\fnms{Denis} \snm{Belomestny}\corref{}\thanksref{t1}\ead[label=e1]{denis.belomestny@uni-due.de}}
\runauthor{D. Belomestny}
\thankstext{t1}{Supported by the Deutsche
Forschungsgemeinschaft through the SPP 1324 ``Mathematical methods for
extracting quantifiable information from complex systems'' and by
Laboratory for Structural Methods of Data Analysis in Predictive
Modeling, MIPT, RF government Grant, ag. 11.G34.31.0073.}
\affiliation{Duisburg-Essen University}
\address[A]{Duisburg-Essen University\\
Altendorferstrasse 11\\
45127 Essen\\
Germany\\
\printead{e1}} 
\end{aug}

\received{\smonth{2} \syear{2012}}
\revised{\smonth{8} \syear{2012}}


\begin{abstract}
In this paper we consider a method of solving optimal stopping problems
in discrete and continuous time based on their dual representation. A
novel and generic simulation-based optimization algorithm not involving
nested simulations is proposed and studied. The algorithm involves the
optimization of a genuinely penalized dual objective functional over a
class of adapted martingales. We prove the convergence of the proposed
algorithm and demonstrate its efficiency for optimal stopping problems
arising in option pricing.
\end{abstract}

%


\end{frontmatter}
%
\section{Introduction}
Let $ (\Omega,\mathcal{F},(\mathcal{F}_{t})_{t\geq0},\P) $ be a
standard filtered probability space, and let $ Z_t $ be an adapted
process satisfying
\[
\E\sup_{t\in[0,T]} |Z_{t}|^{2}<\infty.
\]
Consider the following optimal stopping problem:
%
\begin{equation}
Y^{*}=\sup_{\tau\in\mathcal{T}[0,T]} \E[Z_{\tau}],
\end{equation}
where $ \mathcal{T}[0,T] $ is the set of stopping times taking values
in $ [0,T] $ for some $ T>0. $
Solving the optimal stopping problem \eqref{optstop} is
straightforward in low dimensions. However, many problems arising
in practice have high dimensions, and these applications
have forced the development of simulation-based algorithms for optimal
stopping problems.
There are basically two approaches toward solving \mbox{optimal} stopping
problems: a primal approach and a dual approach. Solving
high-dimensional optimal stopping problems by the primal approach and
Monte Carlo is a challenging task because the determination of the
optimal value function in the primal approach uses a backward dynamic
programming
principle that seems to be incompatible with the forward nature of
Monte Carlo
simulation. Much research was focused on the development of fast
methods to
compute approximations to the optimal value function. One of the most
successful algorithms, and the one adopted most widely by
practitioners, is the Longstaff--Schwartz algorithm. It is based on
approximating the conditional
expectations by the least-squares regression on a given basis of
functions and hence boils down to solving a quadratic optimization problem.
During the last century, the primal approach was, in effect, the only method
available, but in recent years another quite different ``dual''
approach has been discovered by
\citet{Ro} and \citet{HK} that is based on a dual representation for
the optimal value function.
The dual representation involves the minimization of the dual objective
functional over the set of all adapted martingales $ \mathcal{M} $,
where the minimum is attained at some ``optimal'' martingale $ M^{*}
$ that coincides with the martingale in the Doob--Meyer decomposition
of the value process.
In fact, finding such an optimal martingale is as difficult as solving
the original
stopping problem. The so-called martingale duality approach aims at
approximating
the ``optimal'' martingale and then uses this approximation to compute
upper bounds by Monte Carlo.
There are two types of algorithms toward approximating the ``optimal''
martingale $ M^{*}. $ The first one needs a preliminary estimate for
the value process $Y^{*}$ in order to approximate the Doob martingale
$M^{*}.$ The early paper of \citet{AB} uses, for example, the
Longstaff--Schwartz algorithm to construct a pilot estimate for $
Y^{*}$ and then employs sub-simulation to approximate $ M^{*}. $
Another dual algorithm that does not involve sub-simulation, was
suggested in \citet{BBS},
where an approximation for the martingale $ M^{*} $ was constructed
using martingale representation theorem and an approximation of the
value process.
Let us note that the performance of the above two methods deteriorate
sharply as the number of exercise dates increases.
The second type of algorithms is based on the direct optimization of
the dual objective functional over a parameterized set of martingales
and does not require a preliminary estimate of $Y^{*}.$
The recent work of \citet{DFM} uses optimization and sub-simulation to
approximate $ M^{*} $ and $ Y^{*}_{t} $ simultaneously in an
efficient way. However, it becomes less efficient in the case of
continuous optimal stopping problems, as it involves sub-simulations at
each time step. Another ``pure'' dual algorithm was proposed in \citet
{Ro1} and further refined in \citet{Sch}.
Let us finally mention the recent work of \citet{Chr} where a quite
different approach was proposed that uses neither the dual
representation nor Monte Carlo. This approach is based on the excessive
function characterization of the value function for continuous optimal
stopping problems.

The contribution of the current paper is threefold. On the one hand, we
propose a novel dual optimization-based algorithm for solving optimal
stopping problem in discrete and continuous time which does not require
nested Monte Carlo simulations. Our algorithm makes use of the
martingale representation theorem to parametrize the set of martingales
we optimize over. This allows us to obtain surprisingly good results in
a number of benchmark option pricing problems using rather generic sets
of basis functions (trigonometric polynomials) to approximate the
integrand in the martingale representation theorem. In the previous
literature one was able to obtain such bounds only by using either many
sub-simulations or special basis functions [e.g., European deltas in
\citet{BBS} or excessive functions in \citet{Chr}]. On the other hand,
we propose a novel approach toward variance reduction based on the
genuine penalization of the dual objective functional. Last but not the
least, we rigorously analyze the convergence of the proposed dual
algorithm and derive the corresponding convergence rates. Note that as
opposed to the Longstaff--Schwartz algorithm, the convergence of dual
algorithms has not been yet rigorously studied.
Even the convergence of the well-known primal-dual Andersen--Broadie
algorithm is not an obvious issue, as the errors stemming from the
Longstaff--Schwartz algorithm are to be taken into account in a proper
way; see, for example, \citet{B1}.

The paper is organized as follows. In Section~\ref{mainresults} we
formulate the main algorithm and address the convergence issue. In
Section~\ref{MaMaRep} we discuss how to build up a class of martingales
with good approximation properties using the so-called martingale
representation. Section~\ref{sim} contains several numerical examples
illustrating the efficiency of our approach. Section~\ref{concl}
concludes the paper. Finally, in Section~\ref{proofs} the proofs of the
main results together with some auxiliary results are collected. In
particular, we derive a novel concentration inequality for some
empirical process
over parameterized classes of martingales.

\section{Main results}
\label{mainresults}
\subsection{Empirical penalized dual algorithm}
Consider the following optimal stopping problem:
%
\begin{equation}
\label{optstop} Y^{*}_{t}=\mathop{\esssup}_{\tau\in\mathcal{T}[t,T]}
\operatorname {E}[Z_{\tau
}|\mathcal{F}_{t}],\qquad t\in[0,T],
\end{equation}
where $ \mathcal{T}[t,T] $ is the set of stopping times taking values
in $ [t,T] $ for some $ T>0. $
Let $ \mathcal{A} $ stand for the space of all adapted martingales
starting at $0$, then we have the following
dual representation [see \citet{Ro}] for the value process $
Y^{*}_{t}: $
%
\begin{equation}
\label{dualrepr} Y^{*}_{t}=\inf_{M\in\mathcal{A}}
\Bigl(M_{t}+\operatorname{E} \Bigl[ \sup_{u\in[t,T]}(Z_{u}-M_{u})
\big\vert \mathcal{F}_{t} \Bigr] \Bigr).
\end{equation}
The infimum is attained by taking $M=M^{*}, $
where
\[
Y^{*}_{t}=Y^{*}_{0}+M^{*}_{t}-A^{*}_{t}
\]
is the Doob--Meyer decomposition of the supermartingale $ Y_{t}^{*},
$ $ M^{*} $ being a martingale and $ A^{*} $
being an increasing process with $ A^{*}_0 = 0 $. Moreover, the identity
%
\begin{equation}
\label{dualrepras} Y^{*}_{t}= M^{*}_{t}+
\sup_{u\in[t,T]}\bigl(Z_{u}-M^{*}_{u}\bigr)
\end{equation}
holds for all $t\in[0,T]$ with probability $ 1. $
Hence, for an arbitrarily chosen adapted martingale $ M $ with $
M_{0}=0$, the
value
%
\begin{equation}
\label{obj} \operatorname{E} \Bigl[\sup_{u\in[0,T]}(Z_{u}-M_{u})
\Bigr]
\end{equation}
defines an upper bound for $ Y_{0}^{*} $, and the upper bound
will be tight if $ M $ minimizes~\eqref{obj}. On the other hand, we
are interested in martingales $M$ leading to the random variable $
\sup_{t\in[0,T]}(Z_{t}-M_{t})$ with a low variance, since this would
imply faster convergence of a Monte Carlo estimate for \eqref{obj}. By
compromising both requirements, one ends up with the optimization problem
%
\begin{equation}
\label{vardualmin} \inf_{M\in\mathcal{A}} \Bigl\{\E \Bigl[\sup_{t\in
[0,T]}(Z_{t}-M_{t})
\Bigr]+\lambda\sqrt{\Var \Bigl[\sup_{t\in
[0,T]}(Z_{t}-M_{t})
\Bigr]} \Bigr\},
\end{equation}
where $\lambda$ is a nonnegative number determining the degree of
penalization by the variance. Note that due to \eqref{dualrepras} the
Doob martingale $M^{*}$ is one solution of the optimization problem
\eqref{vardualmin}.

Fixing a set of martingales $\mathfrak{M}\subset\mathcal{A}$ and
replacing the true quantities in \eqref{vardualmin} by their
empirical counterparts, we arrive at the following empirical
optimization problem:
%
\begin{equation}
\label{YnPen} M_{n}=\mathop{\arginf}_{M\in\mathfrak{M}} \Biggl( \frac{1}{n}
\sum_{j=1}^{n} Z^{(j)}(M)+\lambda
\sqrt{V_{n}(M)} \Biggr),\qquad \lambda>0,
\end{equation}
where $Z^{(j)}(M),$ $j=1,\ldots, d$ are i.i.d random variables
having the same distribution as
%
\begin{equation}
\label{ZM} Z(M) =\sup_{s\in[0,T]}(Z_{s}-M_{s})
\end{equation}
and
%
\begin{equation}
\label{VM} V_{n}(M)=\frac{1}{n(n-1)}\sum
_{1\leq i<j\leq n} \bigl(Z^{(i)}(M)-Z^{(j)}(M)
\bigr)^{2}.
\end{equation}
The approach based on the empirical optimization problem \eqref{YnPen}
has several obvious advantages. First, it delivers ``true'' upper bound
without use of sub-simulation, thus resulting in a nonnested Monte
Carlo. Second, it does not exclusively focus on finding Doob martingale
and takes advantage of the richness [see \citet{Sch}] of the class $
\mathcal{A}^{*}$ of adapted martingales starting at $0$ and satisfying
%
\begin{equation}
\label{dualrepras0} Y^{*}=\sup_{t\in[0,T]} (Z_{t}-M_{t}),
\qquad\mbox{a.s.}
\end{equation}
Another useful feature of our algorithm which will be proved in the
next section is that the variance
of the r.v. $Z(M_{n})=\sup_{s\in[0,T]}(Z_{s}-M_{n,s})$ is, with high
probability, bounded by a multiple of the r.v.
\[
\inf_{M\in\mathfrak{M}, M'\in\mathcal{A}^{*}}d\bigl(M,M'\bigr),
\]
where $d$ is a deterministic metric on $\mathcal{A}.$ The above
property implies that the variance of $Z(M_{n})$ can be made
arbitrary small by considering classes of martingales $\mathfrak{M}$
with better approximation properties with respect to the solution class
$\mathcal{A}^{*}.$ Last but not least, our approach is applicable to
the case of continuous optimal stopping problems, as it does not
involve regression (or subsimulations) at each discretization step as
in other approaches based on the dynamic programming formulation.
\subsection{Convergence}
Let $ (\Psi,\rho) $ be a metric space. Furthermore, let $ \mathcal
{M}=\{ M(\psi)\dvtx \psi\in\Psi\} $ be a family of adapted continuous
local martingales defined on $ (\Omega, \mathcal{F}, \P). $
%
\begin{defi}
A quadratic $ \rho$-modulus $ \| \mathcal{M} \|_{\rho} $ of a
family $ \mathcal{M}=\{ M(\psi)\dvtx \break \psi\in\Psi\} $ of
continuous local martingales is defined as an $ \mathbb{R}_{+}\cap\{
\infty\} $-valued stochastic process
$ t\mapsto\| \mathcal{M} \|_{\rho,t} $ given by
\[
\| \mathcal{M} \|_{\rho,t}=\mathop{\mathop{\operatorname{ess\,sup}}_{\psi,\phi\in\Psi}}_{
\psi
\neq\phi}
\frac{\sqrt{\langle M(\psi)-M(\phi)\rangle_{t}}}{\rho
(\psi
,\phi)},\qquad t\in[0,T],
\]
where $ \langle M\rangle$ stands for the quadratic variation process
of the continuous local
martingale $ M $.
\end{defi}
For a given subset $ \widetilde\Psi$ of the metric space $ (\Psi
,\rho) $ denote by $ N(\varepsilon,\widetilde\Psi,\rho) $ the
smallest number of closed balls, with $ \rho$-radius $ \varepsilon
>0, $ which cover the set $ \widetilde\Psi$ and define
\[
J(\delta)=\int_{0}^{\delta}\sqrt{\log\bigl[1+ N(
\varepsilon, \widetilde\Psi, \rho)\bigr]} \,d\varepsilon\vadjust{\goodbreak}
\]
for all $\delta>0$.
Denote also by $\mathcal{M}^{*}=\{ M(\psi)\dvtx \psi\in\Psi^{*} \} $ a
subset of $\mathcal{M}$ containing all martingales $M$ that fulfill
\eqref{dualrepras}. In the sequel we shall assume that the family~$
\mathcal{M}$ is rich enough so that
$\mathcal{M}^{*}$ is not empty.
Let us now formulate the main result on the convergence of $ \E
[Z(M_{n})] $ for $M_{n}$ defined in \eqref{YnPen}.
%
\begin{thmm}
\label{conv}
Let $ \mathfrak{M}=\{ M(\psi)\dvtx \psi\in\widetilde\Psi\} $ be a
family of continuous local martingales
satisfying $\| \mathfrak{M} \|_{\rho,T}\leq\Theta$ almost surely,
for some finite $\Theta.$ Let also $\psi^{*}$ be an element of $
\Psi^{*}$ such that $\rho(\psi,\psi^{*})\leq\sigma$ for all $
\psi
\in\widetilde\Psi$ and some $\sigma<\infty.$
Set
\[
\mathfrak{C}=\mathfrak{C}(\widetilde\Psi)=\int_{0}^{\sigma
}
\varepsilon^{-1}J(\varepsilon)\sqrt{\log\bigl[1+ N(\varepsilon,
\widetilde\Psi, \rho )\bigr]} \,d\varepsilon,
\]
and assume that $\mathfrak{C}<\infty.$
Fix some $\varkappa>0$ and $ 0<\delta<1 $ with $ J(1)\log
(1/\delta)\leq\sqrt{n} $, and define
%
\begin{equation}
\label{Mn} M_{n}=\mathop{\arginf}_{M\in\mathfrak{M}} \Biggl( \frac{1}{n}\sum
_{j=1}^{n} Z^{(j)}(M)+\bigl(
\varkappa+\lambda_{n}(\delta/4)\bigr) \sqrt{V_{n}(M)} \Biggr),
\end{equation}
where $Z^{(j)}(M),$ $j=1,\ldots,n,$ and $V_{n}(M)$ are defined in
\eqref{ZM} and \eqref{VM}, respectively,
and $ \lambda_{n}(\alpha)=4(2\sqrt{2\log(2/\alpha)}+\mathfrak
{C})/\sqrt{n} $ for any $ \alpha>0. $
Then it holds for some constant $C>0$ (not depending on $\delta,$
$
n$ and $\varkappa$) with probability at least $1-\delta$,
%
\begin{eqnarray}
\label{biasbound} 0\leq Y(M_{n})-Y^{*}&\leq& C \bigl(
\varkappa+2\lambda_{n}(\delta/4)\bigr) \inf_{\psi\in\widetilde\Psi}\mathcal{R}
\bigl(\psi,\psi^{*}\bigr),
\\
\label{varbound} \sqrt{V(M_{n})}&\leq& C \biggl(1+\frac{2\lambda_{n}(\delta
/4)}{\varkappa
}
\biggr) \inf_{\psi\in\widetilde\Psi}\mathcal{R}\bigl(\psi,\psi^{*}\bigr),
\end{eqnarray}
where $ Y(M)=\E [ \sup_{s\in[0,T]}(Z_{s}-M_{s})  ], $ $
V(M)=\Var [ \sup_{s\in[0,T]}(Z_{s}-M_{s})  ]$
and
\[
\mathcal{R}\bigl(\psi,\psi^{*}\bigr)=\rho\bigl(\psi,\psi^{*}
\bigr)\sqrt{1 \vee\bigl\llvert \log \bigl(\rho\bigl(\psi,\psi^{*}\bigr)
\bigr)\bigr\rrvert }
\]
for any $\psi\in\Psi.$
\end{thmm}
%
\begin{rem}
Note that $Y(M_{n})$ and $V(M_{n})$ are random variables measurable
w.r.t. the $\sigma$-algebra generated by the paths used to compute $
M_{n}.$
\end{rem}
%
\begin{rem}
The condition $\mathfrak{C}<\infty$ roughly means that
$J(\varepsilon)=O(\varepsilon^{1/2+\delta})$ as $\varepsilon\to
0$
for some $\delta>0.$
\end{rem}
%

\textit{Discussion}.
Theorem~\ref{conv} shows that the martingale $ M_{n} $ delivered by
our algorithm has a nice property that the corresponding approximation
error $Y(M_{n})-Y^{*}$ and the square root variance $\sqrt {V(M_{n})}$ can be bounded
from above with high probability by the
quantities proportional to the smallest distance between the classes of
martingales $\mathfrak{M}$ and $\mathcal{A}^{*}$ as measured by $
\rho.$
Hence, if the set $ \mathfrak{M} $ contains at least one martingale
solving \eqref{dualrepras} we get, as expected, $ Y(M_{n})=Y^{*} $
with probability $ 1-\delta. $ In general, the larger is the class $
\mathfrak{M},$ the smaller is the above distance.
However, if the class $ \mathfrak{M} $ is infinite-dimensional,
maximizing the empirical objective functional in \eqref{Mn}
over $ \mathfrak{M} $ may not be well defined or even if $ M_{n} $
exists, it might be difficult to compute.
Instead, one can restrict the maximization to a sequence of
finite-dimensional approximating spaces $ \mathfrak{M}_{n}=\{ M(\psi
)\dvtx \psi\in\Psi_{n} \} $ such that
$ \bigcup_{n} \Psi_{n}$ is dense in $ \Psi^{*}. $ Such a
sequence of approximating spaces is usually called a \textit{sieve}. We
are interested in sieves that are
compact, nondecreasing ($ \mathfrak{M}_{n}\subset\mathfrak
{M}_{n+1}\subset\cdots\subset\mathfrak{M}$) and such that for any
$
n\in\mathbb{N} $ and some $\psi^{*}\in\Psi^{*}$ there exists an
element $ \pi_{n}\psi^{*} $ in $ \Psi_{n} $ satisfying $ \rho
(\psi^{*},\pi_{n}\psi^{*})\to0 $ as $ n\to\infty, $ where $ \pi_{n} $
can be regarded as a projection of $ \psi^{*} $ to $ \Psi_{n}. $
For such sieves Theorem~\ref{conv} implies that
%
\begin{equation}
\label{biasvar} V(M_{n})\stackrel{\P} {\longrightarrow} 0,\qquad n\to\infty,
\end{equation}
provided $ \mathfrak{C}(\Psi_{n})/\sqrt{n}$ remains bounded as $
n\to\infty. $ In the next section we discuss how to get
the martingale sieves $\mathfrak{M}_{n}$ in a constructive way. The
asymptotic relation~\eqref{biasvar} implies that the variance of the
Monte Carlo estimate of $Y(M_{n})$,
\[
Y_{m}(M_{n})=\frac{1}{m}\sum
_{j=n+1}^{n+m} Z^{(j)}(M_{n}),
\]
based on a new, independent sequence of r.v.
\[
\bigl(Z^{(n+1)}(M_{n}),\ldots, Z^{(n+m)}(M_{n})
\bigr)
\]
has the standard deviation of order $o(1/\sqrt{m})$ as $m,n\to
\infty.$ Therefore one can speak about fast convergence rates in this
situation. Let us mention at this place that the primal-dual algorithm of
\citet{AB} has the same variance ``self-reduction'' property [see
\citet
{CG}]: nearer is the preliminary regression estimate $Y_{t}$ of the
value function to the true one, the lower variance has the r.v.
$\sup_{t\in[0,T]}(Z_{t}-M_{t})$ with $M$ based on $Y.$ However, the
results on the speed of the variance decay
in dependence on the number of basis functions and Monte Carlo paths
used in regression step are not yet available in the literature.
%
\begin{rem}
If the class $\widetilde\Psi$ is of Vapnik--Cervonenkis type, that is,
\[
N(\varepsilon, \widetilde\Psi, \rho)\lesssim\varepsilon^{-\beta},\qquad
\varepsilon\to0
\]
for some $\beta>0,$ then the quantity $J(\delta)$ is finite for any
$\delta>0$.
\end{rem}
%
\begin{rem}
\label{nonuniqex}
A natural question is whether the bounds of Theorem~\ref{conv} can be
achieved without using the penalization by empirical variance. The
answer is, in general, no.
To see this, let~$Z_{t}$ be an uniformly\vadjust{\goodbreak} integrable submartingale.
Then~$Z_{t}$ admits the so-called Doob--Meyer decomposition
\[
Z_{t}=Z_{0}+M_{t}+A_{t},
\]
where $M_{t}$ with $M_{0}=0$ is a uniformly integrable martingale,
and $A_{t}$ is an increasing predictable process.
Using the optional sampling theorem, we derive
\[
Y^{*}=\sup_{\tau\in\mathcal{T}[0,T]}\E[Z_{\tau}]=\E [Z_{T}]=Z_{0}+
\E[A_{T}].
\]
Define $M_{t}^{*}=M_{t}+\E[A_{T}|\mathcal{F}_{t}]-\E[A_{T}],$ then
$
Y^{*}=\sup_{t\in[0,T]}(Z_{t}- M^{*}_{t})$
with probability $1.$
Furthermore, the martingale $\widetilde M=M$ fulfills
\[
Y^{*}=\E \Bigl[\sup_{t\in[0,T]}(Z_{t}-\widetilde
M_{t}) \Bigr],
\]
and if $A_{T}$ is not deterministic, then $Y^{*}\neq\sup_{t\in
[0,T]}(Z_{t}- \widetilde M_{t})=Z_{0}+A_{T}$ with positive probability.
Hence, $\widetilde{M}$ solves, along with $M^{*}$, the original
dual problem~\eqref{dualrepr}, but does not have the almost sure
property \eqref{dualrepras}.
Consider now the empirical optimization problem
\[
M_{n}=\mathop{\arginf}_{M\in \{M^{*},\widetilde{M} \}} \Biggl( \frac
{1}{n}\sum
_{j=1}^{n} Z^{(j)}(M) \Biggr)
\]
with $Z(M)=\sup_{t\in[0,T]}(Z_{t}-M_{t}).$ Due to CLT, it obviously holds
\[
\liminf_{n\to\infty}\P(M_{n}=\widetilde M)=\liminf_{n\to\infty
}
\P \Biggl(\sum_{j=1}^{n}
\xi_{j}<0 \Biggr)>0,
\]
where $\xi_{1},\ldots, \xi_{n}$ are i.i.d. random variables
distributed as $A_{T}-\E[A_{T}]. $ Therefore
\[
V(M_{n})=V(\widetilde M)=\Var \Bigl[\sup_{t\in
[0,T]}(Z_{t}-
\widetilde M_{t}) \Bigr]=\Var[A_{T}]>0
\]
with positive probability for any natural number $n$ and the bound
\eqref{varbound} does not hold any longer.
\end{rem}
%
\section{Martingales via martingale representations}
\label{MaMaRep}
Suppose that $ Z_{t}=G_{t}(X_{t}), $ where $ G_{t}\dvtx \mathbb
{R}^{d}\to\mathbb{R} $ is a H\"older function on $[0,T]\times
\mathbb
{R}$ and $ X_{t} $ is a $ d $-dimensional Markov process solving
the following system of SDE's:
%
\begin{equation}
\label{Xsys} dX_t=\mu(t,X_t) \,dt+\sigma(t,X_t)
\,dW_t,\qquad X_0=x.
\end{equation}
The coefficient functions $\mu\dvtx [0,T]\times\mathbb{R}^{d} \rightarrow
\mathbb{R}^d$ and $\sigma\dvtx  [0,T]\times\mathbb{R}^d\rightarrow
\mathbb{R}^{d\times m}$ are supposed to be Lipschitz in space and
$1/2$-H\"older continuous in time, with $m$ denoting the dimension
of the Brownian motion $W=(W^1,\ldots, W^m)^\top$ under
measure $\P$.
It is well known\vadjust{\goodbreak} that under the assumption that a martingale $M_t$ is
square integrable and is adapted to the filtration generated by $
W_{t},$ there
is a square integrable (row vector valued) process
$H_t=(H_{t}^1,\ldots,H_{t}^m)$ satisfying
%
\begin{equation}
\label{MR} M_{t}= \int_0^{t}
H_s \,dW_s.
\end{equation}
It is not hard to see that in the Markovian setting $
Y^{*}_{t}=V(t,X_{t})$, it holds $ H_{s}=\psi(s,X_{s}) $ for some
vector function $ \psi(s,x)=(\psi_{1}(s,x),\ldots,\psi_{m}(s,x))
$ satisfying
\[
\int_{0}^{T}\E\bigl[\bigl|\psi(s,X_{s})\bigr|^{2}
\bigr] \,ds<\infty.
\]
As a result,
\[
M_{t}=M_{t}(\psi)= \int_0^{t}
\psi(s,X_{s}) \,dW_s.
\]
Thus, the set of adapted square-integrable martingales can be
``parameterized'' by the set $ L_{2,P}([0,T]\times\mathbb{R}^{d})$
of square-integrable $m$-dimensional vector functions
$ \psi$ on $ [0,T]\times\mathbb{R}^{d} $ that satisfy $ \| \psi
\|^{2}_{2,\P}:=\int_{0}^{T}\E[|\psi(s,X_{s})|^{2}] \,ds<\infty. $ Let
${\Psi}^{*}$ be a set
of $\psi\in L_{2,P}([0,T]\times\mathbb{R}^{d})$ such that $
M_{t}(\psi)$ solves \eqref{dualrepras}.
Choose a family of finite-dimensional linear models of functions,
called sieves,
with good approximation properties. We consider linear sieves of the form
%
\begin{equation}
\label{SK} {\Psi}_{K}=\{ \beta_{1}\phi_{1}+
\cdots+\beta_{K}\phi_{K}\dvtx \beta_{1},\ldots,
\beta_{K}\in\mathcal{C} \},
\end{equation}
where $ \phi_{1},\ldots,\phi_{K} $ are some given vector functions
with components from the space of bounded continuous functions $
C_{b}([0,T]\times\mathbb{R}^{d})$, and $\mathcal{C}$ is a compact
set in~$\mathbb{R}.$
Next define a class of adapted square-integrable martingales via
\[
\mathfrak{M}_{K}=\bigl\{ M_{t}(\psi)\dvtx \psi\in{
\Psi}_{K}\bigr\}
\]
and set
%
\begin{equation}
\label{Y0nPK} M_{n}:=\mathop{\arginf}_{M\in\mathfrak{M}_{K_{n}}} \Biggl( \frac{1}{n}
\sum_{j=1}^{n}Z^{(j)}(M)+(\varkappa+
\lambda_{n}) \sqrt{V_{n}(M)} \Biggr),
\end{equation}
where $ K_{n}\to\infty$ as $ n\to\infty. $ As can be easily seen
\begin{eqnarray*}
\sqrt{\bigl\langle M-M'\bigr\rangle_{T}} &\leq& \sqrt{T}
\cdot\sup_{x\in
\mathbb
{R}^{d}}\sup_{t\in[0,T]}\bigl|\psi(t,x)-\psi'(t,x)\bigr|
\\
&=&\sqrt{T}\cdot\rho\bigl(\psi,\psi'\bigr)
\end{eqnarray*}
with $ M_{t}=M_{t}(\psi) $ and $ M'_{t}=M_{t}(\psi') $ for any $
\psi,\psi'\in C_{b}([0,T]\times\mathbb{R}^{d})\times\cdots\times
C_{b}([0,T]\times\mathbb{R}^{d}).$ Hence the quadratic $\rho$
-modulus of the family $\mathfrak{M}_{K}$
is bounded by $\sqrt{T}$ with probability $1.$ For many linear
sieves of the form \eqref{SK} and diffusion processes $X,$ it holds
that
\[
\log\bigl[1+N(\varepsilon, {\Psi}_{K}, \rho)\bigr]\lesssim
K^{d+1} \log (1/\varepsilon),\qquad \varepsilon\to0\vadjust{\goodbreak}
\]
and in this situation we have with probability at least $1-\delta$
\[
\sqrt{V(M_{n})}=O (a_{n} ),
\]
where $ a_{n}=\inf_{\psi\in{\Psi}_{K_{n}}, \psi^{*}\in{\Psi}^{*}}
\rho(\psi,\psi^{*}),$ provided
$K_{n}^{d+1}/\sqrt{n}=O(1)$ for $n\to\infty.$

\section{Numerical study}
\label{sim}
In this section we test our algorithm on several benchmark examples
related to American/Bermudan option pricing problems arising in
finance. Let us first give some general details on the implementation
of our algorithm.
First, we need to construct a set of approximating martingales. To this
end we are going to use the martingale representation theorem as
described in Section~\ref{MaMaRep}. It is known [see, e.g., \citet
{BBS}] that in the Markovian setting $Y^{*}_{t}=V(t,X_{t})$ and under
some rather general assumptions on the diffusion process $X$ in
\eqref
{Xsys} the Doob martingale $M^{*}$ with $M^{*}_{0}=0$ has a representation
%
\begin{equation}
\label{DMart} M_{t}^{*}=\int_{0}^{t}
\sum_{i=1}^{d}\frac{\partial
V(u,X_{u})}{\partial
X^{i}}
\sigma^{i}(u,X_{u}) \,dW_{u}.
\end{equation}
Fix now some linear space $\widetilde\Psi$ of functions $\psi\dvtx
[0,T]\times\mathbb{R}^{d}\to\mathbb{R}^{d} .$ The equality~\eqref
{DMart} motivates us to consider the following optimization problem:
%
\begin{eqnarray}
\label{optimynsmooth}&& \psi_{n,\lambda}:=\mathop{\arginf}_{\psi\in\widetilde\Psi} \Biggl\{
\frac
{1}{n}\sum_{j=1}^{n}
Z^{(j)}(\psi)
\nonumber
\\[-8pt]
\\[-8pt]
\nonumber
&&\hspace*{70pt}{}+\lambda\sqrt{\frac{1}{n(n-1)} \sum
_{1\leq i<j\leq n}\bigl(Z^{(i)}(\psi )-Z^{(j)}(\psi)
\bigr)^{2}} \Biggr\}
\end{eqnarray}
with
%
\begin{eqnarray}
\label{Zpsi} Z^{(j)}(\psi)&:=&\sup_{t\in[0,T]} \bigl[G_{t}
\bigl(X^{(j)}_{t}\bigr)-M^{(j)}_{t}(\psi)
\bigr],
\\
\label{Mpsi} M^{(j)}_{t}(\psi)&:=&\int_{0}^{t}
\sum_{i=1}^{d}\sigma^{i}
\bigl(u,X^{(j)}_{u}\bigr)\psi_{i}
\bigl(u,X^{(j)}_{u}\bigr) \,dW^{(j)}_{u}
\end{eqnarray}
and some $\lambda>0,$
where $(W^{(j)}_{t},X^{(j)}_{t})\in\mathbb{R}^{m}\times\mathbb
{R}^{d}, j=1,\ldots,n,$ is the set of trajectories obtained, for
example, by discretizing the system of SDEs \eqref{Xsys}.

\begin{rem}
The construction of the class \eqref{Mpsi} of approximating martingales
is based on some prior information on the underlying process in form of
the matrix $\sigma.$ Moreover, this construction implies that we are
actually aiming at approximating the Doob martingale $M^{*}$ in this case.
\end{rem}
%
\begin{rem}
Let us discuss the choice of the penalization parameter $\lambda$ in
more details. On the one side,
the parameter $\lambda$ can be chosen according to Theorem~\ref
{conv}, that is, $\lambda=\varkappa+4(2\sqrt{2\log(2/\alpha
)}+\mathfrak
{C})/\sqrt{n}$ for some $\varkappa>0$ and $\alpha>0.$ This choice,
however, requires knowledge of
\[
\mathfrak{C}=\mathfrak{C}(\widetilde\Psi)=\int_{0}^{\sigma
}
\varepsilon^{-1}J(\varepsilon)\sqrt{\log\bigl[1+ N(\varepsilon,
\widetilde\Psi, \rho )\bigr]} \,d\varepsilon,
\]
which might be difficult to compute in concrete situations. On the
other side, $\lambda$ can be found empirically by minimizing the
``out of sample'' variance and mean of the r.v. $Z(M_{n}).$
This would require some additional computational efforts.
\end{rem}
In all examples below we use the Euler scheme and $n_{\mathrm{disc}}=200$
discretization points to approximate \eqref{Xsys}. The integral in
\eqref{Mpsi} can be then easily approximated through the sum
\[
\sum_{l=1}^{n_{\mathrm{disc}}}\sum
_{i=1}^{d}\sigma^{i}\bigl(u_{l},X^{(j)}_{u_{l}}
\bigr)\psi_{i}\bigl(u_{l},X^{(j)}_{u_{l}}
\bigr) \bigl(W^{(j)}_{u_{l+1}}-W^{(j)}_{u_{l}}
\bigr).
\]
As to the choice of linear space $\widetilde\Psi,$ we are striving
for the most generic choice not involving special functions like
European deltas as in \citet{BBS}. In all examples to follow we first
make a basic variable transformation and
then use trigonometric bases. Let us also comment on the optimization
problem \eqref{optimynsmooth}
which is convex (at least for $n$ large enough), provided $
\widetilde
\Psi$ is a linear space.
Note, however, that the objective functional in \eqref{optimynsmooth}
is, in general, not smooth. In order to avoid computational problems
related to the nonsmoothness of $ Z(\psi), $ we smooth it [see
\citet
{N} for some theoretical justification] and consider instead $Z$ the
functional
%
\begin{equation}
\label{Zppsi} Z_{p}(\psi)=p^{-1}\log \biggl( \int
_{0}^{T}\exp \bigl(p\bigl(G_{s}(X_{s})-M_{s}(
\psi )\bigr) \bigr) \,ds \biggr),
\end{equation}
where $M_{t}(\psi)=\int_{0}^{t}\sum_{i=1}^{d}\sigma^{i}(u,X_{u})\psi_{i}(u,X_{u}) \,dW^{i}_{u}.$
An alternative expression for $ Z_{p}(\psi) $ is
%
\begin{equation}
\label{zp} Z_{p}(\psi)=Z(\psi)+p^{-1}\log \biggl( \int
_{0}^{T}\exp \bigl(p\bigl(Z_{s}-M_{s}(
\psi)-Z(\psi)\bigr)\bigr) \,ds \biggr).
\end{equation}
It follows from representation \eqref{zp} that
\[
0\leq Z_{p}(\psi)-Z(\psi)\leq p^{-1}\log T.
\]
Hence $ Z_{p}(\psi)\to Z(\psi) $ as $ p\to\infty. $ The advantage
of using $ Z_{p}(\psi) $ instead of $ Z(\psi) $ is that the
standard gradient-based optimization routines can be used to compute $
\psi_{n,\lambda}$.

\subsection{American put on a single asset}
\label{S1DPut}
We start with analyzing the continuously exercisable American put
option on a single asset, the simplest American-type option. We assume
the asset price follows the geometric Brownian motion process
\[
dX_{t}=rX_{t} \,dt+\sigma X_{t}
\,dW_{t},
\]
where $r=0.06,$ $\sigma=0.4,$ $W_{t}$ is the standard Brownian
motion, and the stock pays no dividends. The option has a strike price
of $K=100$ and a maturity of $T=0.5,$ and the payoff upon
exercise at time $t$ is $G(X_{t})=e^{-rt}(K-X_{t})^{+}.$

In our implementation of \eqref{optimynsmooth} we take $\widetilde
\Psi_{L}$ to be a linear space of functions $\psi\dvtx  [0,T]\times
\mathbb
{R}^{d}\mapsto\mathbb{R} $ such that
\[
\psi(t,x)\in\operatorname{span} \bigl\{\zeta_{k}\bigl(y_{t}(x)
\bigr),\xi_{k}\bigl(y_{t}(x)\bigr), k=0,\ldots,L \bigr\},
\]
where $y_{t}(x)=\frac{1}{T-t}\log(x/K)$ and
\begin{eqnarray*}
\zeta_{k}(z)&=& \cases{ 0, &\quad $z<-0.5,$ \vspace*{2pt}
\cr
\sin(k\cdot z),
&\quad $|z|\leq0.5,$ \vspace*{2pt}
\cr
1, & \quad$z>0.5,$}
\\
\xi_{k}(z)&=& \cases{ 0, & \quad $z<-0.5,$ \vspace*{2pt}
\cr
\cos(k\cdot z), &\quad
$|z|\leq0.5,$ \vspace*{2pt}
\cr
1, &\quad  $z>0.5.$}
\end{eqnarray*}
%
\begin{table}
\caption{Upper bounds for the standard one-dimensional American put with
parameters $K=100, r=0.06, T=0.5$ and $\sigma=0.4$ obtained using the
linear space $\widetilde\Psi_{5}$. Values for two different values
of the parameter $\lambda$ are presented}\label{1Dput}
\begin{tabular*}{\textwidth}{@{\extracolsep{\fill}}ld{2.4}d{2.14}d{2.14}c@{}}
\hline
\multicolumn{1}{@{}l}{$\bolds{X_{0}}$} &\multicolumn{1}{c}{\textbf{True value}} &
\multicolumn{1}{c}{\textbf{Upper bound} $\bolds{Y_{10^{4},0}(10^{5})}$} &
\multicolumn{1}{c}{\textbf{Upper bound} $\bolds{Y_{10^{4},2}(10^{5})}$} & \multicolumn{1}{c@{}}{\textbf{Time (sec)}}
\\
\hline
\phantom{0}80 & 21.6059 & 21.63044\ (0.04354) & 21.64156\ (0.01321) & 53 \\
\phantom{0}90 & 14.9187 & 14.92159\ (0.01750) & 14.93001\ (0.00576) & 51\\
100 & 9.9458 & 9.93455\ (0.01354) & 9.94712\ (0.00423) & 47\\
110 & 6.4352 & 6.41561\ (0.01329) & 6.42911\ (0.00479) & 47\\
120 & 4.0611 & 4.03417\ (0.01127) & 4.04883\ (0.00392) & 43\\
\hline
\end{tabular*}
\end{table}

Table~\ref{1Dput} is obtained using the following two-step procedure.
First, we generate $n=10\mbox{,}000$ ``training'' paths on which we solve
optimization \eqref{optimynsmooth} to get $\psi_{n,\lambda}$. In
the second step we use $N=100\mbox{,}000$ new paths to test the martingale
resulting from $\psi_{n,\lambda}$ and to get the final estimate
\[
Y_{n,\lambda}(N)=\frac{1}{N}\sum_{j=n+1}^{n+N}
Z^{(j)}(\psi_{n,\lambda}).
\]
The values in Table~\ref{1Dput} are reported together with the standard
deviations obtained by repeating the ``testing'' step $100$ times.
The times in the last column of the table give the duration of the
``training'' step. By inspecting Table~\ref{1Dput} one can draw several
conclusions. First, the values of the upper bound $
Y_{10^{4},2}(10^{5})$ are almost exact. Second, the penalization with
the empirical variance ($\lambda=2$) reduces the standard deviation
by a factor of three. Finally, our approach is able to compete with the
very powerful method of \citet{Chr} (perhaps with a little bit longer
computational time).

\subsection{American puts on the cheapest of $d$ assets}

In this section, we study the performance of our approach for
multiasset American options, where traditional lattice techniques
usually suffer from serious numerical constraints. Specifically, we
price the American put option on the cheapest of $d$ assets. This
example was also studied by \citet{Ro}. The risk-neutral dynamics for
$
d$-dimensional underlying process $X$ is given by
\[
dX_{t}^{i}=rX_{t}^{i} \,dt+
\sigma_{i} X_{t}^{i} \,dW^{i}_{t},\qquad
i=1,\ldots,d,
\]
where $W^{1}_{t},\ldots,W^{d}_{t}$ are $d$ independent Brownian
motions. The payoff at time $t$ is equal to
\[
G(X_{t})=e^{-rt} \Bigl(K-\min_{k=1,\ldots,d}
X^{k}_{t} \Bigr)^{+}.
\]
In our numerical experiment we take $d=2,$ $\sigma_{i}=\sigma
=0.4,$
$r=0.06$ and $K=100$
and consider linear space $\widetilde\Psi_{L}$ of functions $\psi\dvtx
[0,T]\times\mathbb{R}^{d}\mapsto\mathbb{R}^{2} $ such that
\begin{eqnarray}\label{psi1}
\psi_{1}(t,x)&\in&\operatorname{span} \bigl\{
\zeta_{k}\bigl(y_{t}^{1}(x)\bigr),
\xi_{k}\bigl(y_{t}^{1}(x)\bigr),
\zeta_{k}\bigl(y_{t}^{1}(x)\bigr)1
\bigl(y_{t}^{1}\leq y_{t}^{2}\bigr),\nonumber\\
&&\hspace*{24pt}{}\xi_{k}\bigl(y_{t}^{1}(x)\bigr)1
\bigl(y_{t}^{1}\leq y_{t}^{2}\bigr),
\\
 &&\hspace*{24pt}{} \zeta_{k}\bigl(y_{t}^{1}(x)+y_{t}^{2}(x)
\bigr),\xi_{k}\bigl(y_{t}^{1}(x)+y_{t}^{2}(x)
\bigr), k=0,\ldots,L \bigr\}\nonumber
\end{eqnarray}
and
%
\begin{eqnarray}\label{psi2}
\psi_{2}(t,x)&\in&\operatorname{span} \bigl\{
\zeta_{k}\bigl(y_{t}^{2}(x)\bigr),
\xi_{k}\bigl(y_{t}^{2}(x)\bigr),
\zeta_{k}\bigl(y_{t}^{2}(x)\bigr)1
\bigl(y_{t}^{2}\leq y_{t}^{1}\bigr),\nonumber\\
&&\hspace*{26pt}\xi_{k}\bigl(y_{t}^{2}(x)\bigr)1
\bigl(y_{t}^{2}\leq y_{t}^{1}\bigr),
\\
 &&\hspace*{26pt} \zeta_{k}\bigl(y_{t}^{1}(x)+y_{t}^{2}(x)
\bigr),\xi_{k}\bigl(y_{t}^{1}(x)+y_{t}^{2}(x)
\bigr), k=0,\ldots,L \bigr\}\nonumber
\end{eqnarray}
with $\zeta_{k},$ $\xi_{k}$ defined in Section~\eqref{S1DPut}
and $
y^{1}_{t}(x)=\frac{1}{T-t}\log(x^{1}/K),$ $y^{2}_{t}(x)=\frac
{1}{T-t}\log(x^{2}/K).$

\begin{table}
\tabcolsep=0pt
\caption{Upper bounds (with standard deviations) for the $ 2 $-dimensional
Bermudan min-puts with parameters
$K=100, r=0.06$, $\sigma=0.4$}\label{2DMinPut}
\begin{tabular*}{\textwidth}{@{\extracolsep{\fill}}lccccc@{}}
\hline
\multicolumn{1}{@{}l}{$\bolds{X^{1}_{0}}$} &
\multicolumn{1}{c}{$ \bolds{X^{2}_{0}}$} &
\textbf{True value (FD)}& \multicolumn{1}{c}{\textbf{Upper bound} $\bolds{Y_{10^{4},0}(10^{5})}$}&
\multicolumn{1}{c}{\textbf{Upper bound} $\bolds{Y_{10^{4},2}(10^{5})}$}&
\multicolumn{1}{c@{}}{\textbf{Times (sec)}}\\
\hline
\phantom{0}80 &\phantom{0}80 &37.30 &37.65877 (0.02832) & 37.65921 (0.00912) & 67\\
100 &100 &25.06 & 25.16745 (0.02341) & 25.17551 (0.00778) & 63\\
120 &120 &15.92 &15.93370 (0.01949) & 15.94191 (0.00611)& 61\\
\hline
\end{tabular*}
\end{table}

Table~\ref{2DMinPut} is again obtained using a two-step procedure as
described in Section~\ref{S1DPut} and the linear space $\widetilde
\Psi_{7}.$ The results can be significantly improved by adding to $
\widetilde\Psi_{7}$ some special functions, like European deltas or
harmonic functions.

%
\begin{table}[b]
\tabcolsep=0pt
\caption{Bounds (with standard deviations) for 2-dimensional Bermudan max
call with parameters $\kappa=100, r=0.05$, $\sigma=0.2$, $\delta=0.1$}\label{Tab1}
\begin{tabular*}{\textwidth}{@{\extracolsep{\fill}}lccccc@{}}
\hline
\multicolumn{1}{c}{$\bolds{X^{1}_{0}}$} & \multicolumn{1}{c}{$\bolds{X^{2}_{0}}$} & \multicolumn{1}{c}{\textbf{Upper bound} $\bolds{Y_{10^{4},0}(10^{5})}$} &
\multicolumn{1}{c}{\textbf{Upper bound} $\bolds{Y_{10^{4},2}(10^{5})}$}& \multicolumn{1}{c}{\textbf{A\&B Price interval}}  &
\multicolumn{1}{c}{\textbf{Time (sec)}}\\
\hline
\phantom{0}90 & \phantom{0}90 & \phantom{0}8.07742 (0.00832) & \phantom{0}8.08012 (0.00313) & [8.053, 8.082] & 58\\
100 & 100 & 14.01900 (0.01405) & 14.02131 (0.00466) & [13.892, 13.934]
& 61\\
110 & 110 & 21.60967 (0.01798) & 21.62144 (0.00521) & [21.316, 21.359]
& 64\\
\hline
\end{tabular*}
\end{table}
%

\subsection{Bermudan max-calls on $d $ assets}

This is a benchmark example studied in \citet{BG}, \citet{HK} and \citet{Ro}
among others. Specifically, the model with $d$ identically distributed
assets is considered, where each underlying has dividend yield $\delta$.
The risk-neutral dynamic of assets is given by
\[
\frac{dX_{t}^{k}}{X_{t}^{k}}=(r-\delta)\,dt+\sigma \,dW_{t}^{k},\qquad
k=1,\ldots,d,
\]
where $W_{t}^{k}, k=1,\ldots,d$, are independent one-dimensional Brownian
motions and $r,\delta,\sigma$ are constants. At any time $t\in
\{t_{0},\ldots,t_{\mathcal{I}}\}$ the holder of the option may exercise
it and
receive the payoff
\[
G(X_{t})=e^{-rt}\bigl(\max\bigl(X_{t}^{1},\ldots,X_{t}^{d}
\bigr)-K\bigr)^{+}.
\]
We consider a two-dimensional example where $t_{i}=iT/\mathcal{I},
i=0,\ldots,\mathcal{I}$, with $T=3, \mathcal{I}%
=9.$ In order to construct the linear space $\widetilde\Psi_{L}$ we
again use the functions $\psi\dvtx  [0,T]\times\mathbb{R}^{d}\mapsto
\mathbb
{R}^{2} $ with coordinate functions defined in \eqref{psi1} and
\eqref
{psi2}, respectively.
Table~\ref{Tab1} is obtained by setting $L=7.$
One can observe that the results of Table~\ref{Tab1} are especially
good for small values of $X_{0}.$ For example, the upper bound $
Y_{10^{4},2}(10^{5})$ for $X_{0}=(90,90)$ almost coincides with the
exact value $Y_{0}^{*}$ and was previously obtained only by using
either European deltas [see \citet{BBS}] or many sub-simulations; see
\citet{AB}.
As can be seen from Table~\ref{Tab2}, the upper bound [$
X_{0}=(90,\ldots,90)$] remains tight as the dimension $d$ increases.

\begin{table}
\caption{Bounds (with standard deviations) for $d$-dimensional
Bermudan max-call with parameters $\kappa=100, r=0.05$, $\sigma=0.2$,
$\delta=0.1$}%
\label{Tab2}
\begin{tabular*}{\textwidth}{@{\extracolsep{\fill}}lcccc@{}}
\hline
$\bolds{d}$ & \multicolumn{1}{c}{\textbf{Upper bound} $\bolds{Y_{10^{4},0}(10^{5})}$}&
\textbf{Upper bound} & \textbf{A\&B price interval}& \textbf{Time (sec})\\
\hline
3 & 11.28986 (0.00939) & 11.29100 (0.00326) & [11.265, 11.308] & 73\\
5 & 16.68231 (0.01405) & 16.69506 (0.00467) & [16.602, 16.655] & 80\\
\hline
\end{tabular*}
\end{table}

\section{Conclusion}
\label{concl}
This paper proposes an efficient and self-contained dual algorithm for
solving optimal stopping problems in discrete and continuous time which
is based on the direct minimization of the penalized dual objective
functional over a genuinely parameterized set of martingales. We
analyze the asymptotic properties of the estimated value function and
show that its variance can be made arbitrarily small by a proper choice
of approximating martingales. From the methodological point of view,
the probabilistic tools developed in the paper can be used to analyze
the convergence of various types of empirical optimization problems
arising in computational stochastics and finance.

\section{Proofs of main results}
\label{proofs}
\subsection{\texorpdfstring{Proof of Theorem~\protect\ref{conv}}{Proof of Theorem 2.2}}
Let us first sketch the main steps of the proof.
Our main interest lies in estimating the quantities $Y(M_{n})-Y^{*}$
and $V(M_{n}).$ In order to obtain these estimates
we need a kind of uniform (over $M\in\mathfrak{M}$) concentration
inequality for the empirical process
\[
\mathcal{E}_{n}(M)=\frac{1}{n}\sum_{j=1}^{n}
\bigl(Z^{(j)}(M)-\E \bigl[Z(M)\bigr]\bigr)=\frac{1}{n}\sum
_{j=1}^{n} Z^{(j)}(M) -Y(M)
\]
that gives probabilistic bounds for $\sqrt{n}\cdot\mathcal
{E}_{n}(M)$
in terms of the empirical variance $ V_{n}(M). $ Indeed, such an
inequality would allow us to get an upper bound for the quantity $
Y(M_{n})+\varkappa\sqrt{V_{n}(M_{n})}$ with $\varkappa>0$ in terms
of $\mathcal{Q}_{n}(M_{n}),$ where
\[
\mathcal{Q}_{n}(M)=\frac{1}{n}\sum_{j=1}^{n}
Z^{(j)}(M)+\bigl(\varkappa +\lambda_{n}(\delta/2)\bigr)
\sqrt{V_{n}(M)}.
\]
Unfortunately, the usual concentration inequalities could not be used
here, as they would provide us with the bounds in terms of the true
variance $ V(M) $ and not in terms of the empirical one $ V_{n}(M).
$ However, there is another, less-known type of concentration
inequalities for self-normalized empirical processes [see \citet{BGR}],
and this is exactly what we need. We extend the above inequalities to
the case of general family of random variables. As a next step, in
order to derive a bound for $V(M_{n}),$ we need a kind of uniform\vadjust{\goodbreak}
concentration inequality for the empirical process $\Delta_{n}(\psi
)=(V(M(\psi))-V_{n}(M(\psi)))$ that holds uniformly over the set $
\widetilde\Psi$ and gives probabilistic bounds for $\sqrt{n}\cdot
\Delta_{n}(\psi)$ in terms of $\rho(\psi,\psi^{*})$ for any
fixed $
\psi^{*}\in\Psi^{*}.$
The latter type of inequality cannot be derived from the well-known
concentration inequalities for selfbounding random variables [see,
e.g., \citet{DL}], since variance $V(M)$ is a highly nonlinear
function of $M$ and the random variable $Z(M)$ is usually not
bounded. The corresponding concentration inequality making use of the
local subgaussianity of $V(M), $ is presented in Section~\ref{proofs}
and can be interesting in its own right. Finally, using the inequality
$\mathcal{Q}_{n}(M_{n})\leq\mathcal{Q}_{n}(M), $ that holds for any
$M\in\mathfrak{M},$ we will arrive at \eqref{biasbound} and \eqref
{varbound}.

\textit{Part}~1: The following proposition allows us to derive
uniform bounds for the empirical process $\sqrt{n}\cdot\mathcal
{E}_{n}(M)$ in terms of the empirical variance $V_{n}(M).$

\begin{prop}
\label{selfbound}
Let $ \mathfrak{X} $ be a family of centered and normalized random
variables on a common probability space $(\Omega,\mathcal{F},\P)$
with finite bracketing number in $L_{2}(\P)$
such that
\[
\limsup_{n\to\infty}\E \Bigl[ \sup_{X\in\mathfrak{X}}\max\bigl\llvert \sqrt {n}
\cdot\mathbb{E}_{n}[X]\bigr\rrvert \Bigr]\leq\mathfrak{C}<\infty
\]
for some positive constant $ \mathfrak{C}=\mathfrak{C}(\mathfrak{X}),
$ where
\[
\mathbb{E}_{n}[X]= \frac{1}{n%
}\sum
_{j=1}^{n}X^{(j)}
\]
and $X^{(1)},\ldots, X^{(n)}$ are i.i.d. copies of the element $
X\in
\mathfrak{X}.$
Define
\[
W_{n}(X)=\frac{\mathbb{E}_{n}[X]}{\sqrt{V_{n}(X)}}
\]
with
\[
V_{n}(X)=\frac{1}{n}\sum_{j=1}^{n}
\bigl(X^{(j)}\bigr)^{2}.
\]

Then for any $\kappa>0$ and $\alpha>\sqrt{2},$ one can find some
positive $ \theta$
and $ n_{0} $ depending on $ \mathfrak{X}, $ $\alpha$ and $
\kappa
$ such that, for $ n\geq n_{0} $ and for any $ x\in
[0,\theta\sqrt{n}]$
\[
\P \Bigl(\sup_{X\in\mathfrak{X}} \bigl|\sqrt{n}\cdot W_{n}(X)\bigr|\geq(x+ \alpha
\mathfrak{C}) \Bigr)\leq2 \exp \biggl( -\frac{x^{2}}{4\alpha^{2}(1+\kappa)} \biggr).
\]
\end{prop}
For the case of noncentered and nonnormalized random variables $X$,
one can derive from Proposition~\ref{selfbound} the following corollary.
%
\begin{cor}
\label{selfboundcor}
Let $ \mathfrak{X} $ be a family class of random variables on a
common probability space $(\Omega,\mathcal{F},\P)$ with finite
bracketing number in $L_{2}(\P)$
such that
\[
\sup_{X\in\mathfrak{X}}\E\llvert X\rrvert^{2}<\infty\vadjust{\goodbreak}
\]
and
\[
\limsup_{n\to\infty}\E \Bigl[ \sup_{X\in\mathfrak{X}}\bigl\llvert \sqrt{n}\cdot
\mathbb{E}_{n}[X-\E X]\bigr\rrvert \Bigr]\leq\mathfrak{C}<\infty
\]
for some positive constant $ \mathfrak{C}=\mathfrak{C}(\mathfrak{X}).
$
Define
\[
W_{n}(X)=\frac{\mathbb{E}_{n}[X]-\E[X]}{\sqrt{V_{n}(X)}}
\]
with
\[
V_{n}(X)=\frac{1}{n}\sum_{j=1}^{n}
\bigl(X^{(j)}-\mathbb{E}_{n}[X]\bigr)^{2}.
\]

Then for any $\kappa>0$ and $\alpha>\sqrt{2},$ one can find some
positive $ \theta$
and $ n_{0} $ depending on $ \mathfrak{X}, $ $\alpha$ and $
\kappa
$ such that, for $ n\geq n_{0} $ and for any $ x\in
[0,\theta\sqrt{n}]$,
\[
\P \biggl(\sup_{X\in\mathfrak{X}} \bigl|\sqrt{n}\cdot W_{n}(X)\bigr|\geq
\frac
{\sqrt{2}(x+ \alpha\mathfrak{C})}{1-\sqrt{2}(x+ \alpha\mathfrak
{C})/n} \biggr)\leq2 \exp \biggl( -\frac{x^{2}}{4\alpha^{2}(1+\kappa)} \biggr),
\]
provided $\sqrt{2}(x+ \alpha\mathfrak{C})<n.$ As a result, by fixing
some $\delta>0$ with $\log(1/\delta)\leq\sqrt{n}$ and taking $
x=2\alpha\sqrt{(1+\kappa)\log(4/\delta)}$, we get with
probability at
least $1-\delta$
\[
\sup_{X\in\mathfrak{X}} \bigl|\sqrt{n}\cdot W_{n}(X)\bigr|\geq2\sqrt {2}\alpha
\cdot\bigl(2\sqrt{(1+\kappa)\log(2/\delta )}+\mathfrak{C}\bigr)
\]
for all $n>n_{0}.$
\end{cor}

\textit{Part}~2: Next we need the concentration inequality for the
empirical process $\sqrt{n}\cdot(V_{n}(M)-V(M))$. The following
proposition is proved in Section~\ref{expineqfixedpsi}.
%
\begin{prop}
\label{EVexpbound}
Let $ \mathfrak{M}=\{ M(\psi)\dvtx \psi\in\widetilde\Psi\} $ be a
family of continuous local martingales, where $\widetilde\Psi$ is a
subspace of the metric space $(\Psi,\rho).$ Suppose that $\|
\mathfrak{M} \|_{\rho,T}\leq\Theta$ a.s.
for some finite $\Theta$ and
\[
J=\int_{0}^{1}\sqrt{\log\bigl[1+N(\varepsilon,
\widetilde\Psi, \rho)\bigr]} \,d\varepsilon<\infty.
\]
Denote $\Delta_{n}(\psi)=V_{n}(M(\psi))-V(M(\psi))$ for any $
\psi\in
\Psi,$ then
for any fixed $\psi^{*}\in\mathcal{M}^{*} $ such that $ \sup_{\psi
\in\widetilde\Psi}\rho(\psi,\psi^{*})< \infty$ it holds
\[
\P \biggl(\sup_{\psi\in\widetilde\Psi}\biggl\llvert \frac{\sqrt{n}\cdot
\Delta_{n}(\psi)}{\mathcal{R}^{2}(\psi,\psi^{*})}\biggr\rrvert >U
\biggr)\leq \exp \biggl(-\frac{D\cdot U}{J} \biggr)
\]
for any $U>0$ and some constant $D>0$ depending on $\Theta,$ where
\[
\mathcal{R}\bigl(\psi,\psi'\bigr)=\rho\bigl(\psi,\psi'
\bigr)\sqrt{1\vee\bigl\llvert \log \bigl(\rho\bigl(\psi ,\psi'\bigr)
\bigr)\bigr\rrvert }
\]
for any $\psi,\psi'\in\Psi.$
\end{prop}

\textit{Part}~3: Now we can begin with the proof of Theorem~\ref{conv}.
By Corollary~\ref{selfboundcor} it holds for any $\psi\in\widetilde
\Psi$ with probability at least $ 1-\delta/2 $,
\begin{eqnarray*}
Y(M_{n})+\varkappa\sqrt{V_{n}(M_{n})} & \leq&
\frac{1}{n}\sum_{j=1}^{n}Z^{(j)}(M_{n})+
\bigl(\varkappa+\lambda_{n}(\delta/4)\bigr) \sqrt {V_{n}(M_{n})}
\\
& \leq&\frac{1}{n}\sum_{j=1}^{n}Z^{(j)}
\bigl(M(\psi)\bigr)+\bigl(\varkappa+\lambda_{n}(\delta/4)\bigr)
\sqrt{V_{n}\bigl(M(\psi)\bigr)}
\\
& \leq& Y\bigl(M(\psi)\bigr)+\bigl(\varkappa+2\lambda_{n}(\delta/4)
\bigr) \sqrt {V_{n}\bigl(M(\psi)\bigr)}.
\end{eqnarray*}
Proposition~\ref{EVexpbound} implies that with probability at least
$ 1-\delta/4 $,
\begin{eqnarray*}
V_{n}\bigl(M(\psi)\bigr)&\leq& V\bigl(M(\psi)\bigr)+JD^{-1}
\log(4/\delta) \frac
{\mathcal
{R}^{2}(\psi,\psi^{*})}{\sqrt{n}}
\\
&\leq& V\bigl(M(\psi)\bigr)+C \mathcal{R}^{2}\bigl(\psi,
\psi^{*}\bigr)
\end{eqnarray*}
for some universal constant $C$, provided $J\log(1/\delta)\leq
\sqrt {n}.$
Hence, using the elementary inequality $ \sqrt{a+b}\leq\sqrt {a}+\sqrt {b}, $ we get
\[
Y(M_{n})+\varkappa\sqrt{V_{n}(M_{n})} \leq Y
\bigl(M(\psi)\bigr)+\bigl(\varkappa +2\lambda_{n}(\delta/4)\bigr) \bigl[
\sqrt {V\bigl(M(\psi)\bigr)} + \sqrt{C} \mathcal{R}\bigl(\psi,\psi^{*}
\bigr) \bigr]
\]
with probability at least $ 1-3\delta/4. $ By the
Burkholder--Davis--Gundy inequality,
\[
Y\bigl(M(\psi)\bigr)-Y^{*}\leq\Theta\rho\bigl(\psi,\psi^{*}
\bigr)
\]
and $V(M(\psi))\leq\Theta^{2}\rho^{2}(\psi,\psi^{*}) $ for any
$
\psi\in\widetilde\Psi. $ Therefore
\[
Y(M_{n})-Y^{*}\leq2\sqrt{C} \Theta\bigl(1+\varkappa+2
\lambda_{n}(\delta /4)\bigr)\mathcal{R}\bigl(\psi,\psi^{*}
\bigr)
\]
and
\[
\sqrt{V_{n}(M_{n})}\leq2\sqrt{C} \Theta\varkappa^{-1}
\bigl(1+\varkappa +2\lambda_{n}(\delta/4)\bigr)\mathcal{R}\bigl(\psi,
\psi^{*}\bigr).
\]
Using again Proposition~\ref{EVexpbound}, we get with probability at
least $1-\delta$,
\begin{eqnarray*}
\sqrt{V(M_{n})}&\leq&\sqrt{V_{n}(M_{n})}+\sqrt{C}
\mathcal{R}\bigl(\psi ,\psi^{*}\bigr)\\
&\leq&3\sqrt{C} \Theta
\varkappa^{-1}\bigl(1+\varkappa+2\lambda_{n}(\delta/4)\bigr)
\mathcal{R}\bigl(\psi,\psi^{*}\bigr).
\end{eqnarray*}

\textit{Part}~4: To finish the proof of Theorem~\ref{conv}, it
suffices to prove the following proposition.
%
\begin{prop}
\label{donsker}
Let $ \widetilde\Psi$ be a subspace of the metric space $(\Psi
,\rho
)$
such that
$ \rho(\psi,\psi^{*})\leq\sigma$ for some $\psi^{*}\in\Psi^{*},$
all $ \psi\in\widetilde\Psi$ and some $ \sigma>0. $ Define $
\mathfrak{M}=\{ M(\psi)\dvtx \psi\in\widetilde\Psi\} $ and
set
\[
\mathfrak{C}=\int_{0}^{\sigma}\varepsilon^{-1}J(
\varepsilon)\sqrt {\log\bigl[ 1+N(\varepsilon, \widetilde\Psi, \rho)\bigr]} \,d
\varepsilon.
\]
If $\| \mathfrak{M} \|_{\rho,T}\leq\Theta$ a.s. and $ \mathfrak
{C}<\infty, $ then there is a constant $A$ depending on $\Theta,$
such that
\[
\limsup_{n\to\infty}\E \Bigl[ \sup_{M\in\mathfrak{M}}\bigl\llvert \mathbb
{G}_{n}\bigl[Z(M)\bigr]\bigr\rrvert \Bigr]\leq A\mathfrak{C}
\]
with
\[
\mathbb{G}_{n}\bigl[Z(M)\bigr]=\frac{1}{\sqrt{n}}\sum
_{j=1}^{n}\bigl(Z^{(j)}(M)-\E\bigl[Z(M)\bigr]
\bigr).
\]
\end{prop}

\begin{pf}
We follow the proof of Lemma 19.34 in \citet{VV} with some
straightforward modifications. It holds
\begin{eqnarray*}
\mathbb{G}_{n}\bigl(M(\psi)\bigr) & =&\frac{1}{\sqrt{n}}\sum
_{j=1}^{n} \bigl(Z^{(j)}\bigl(M(\psi)\bigr)-
\E\bigl[ Z^{(j)}\bigl(M(\psi)\bigr)\bigr] \bigr)
\\
&=&\frac{1}{\sqrt{n}}\sum_{j=1}^{n}
\bigl(Z^{(j)}\bigl(M(\psi)\bigr)-Z^{(j)}%
\bigl(M\bigl(
\psi^{\ast}\bigr)\bigr) \bigr)
\\
&& {}+\frac{1}{\sqrt{n}}\sum_{j=1}^{n}
\bigl(\E\bigl[ Z^{(j)}\bigl(M\bigl(\psi^{\ast}\bigr)\bigr)\bigr] -
\E\bigl[ Z^{(j)}\bigl(M(\psi)\bigr)\bigr] \bigr),
\end{eqnarray*}
since $\Var[Z(M(\psi^{\ast}))]=0.$ Setting
\[
K_{T}=\sup_{\psi\in
\widetilde{\Psi}}\sup_{t\in\lbrack0,T]}\bigl\llvert
M_{t}(\psi )-M_{t}\bigl(\psi^{*}\bigr)\bigr
\rrvert ,
\]
we derive
\[
\bigl\llvert Z\bigl(M(\psi)\bigr)-Z\bigl(M\bigl(\psi^{\ast}\bigr)\bigr)
\bigr\rrvert \leq K_{T},\qquad \bigl\llvert \E\bigl[Z\bigl(M(\psi)\bigr)
\bigr]-\E\bigl[Z\bigl(M\bigl(\psi^{\ast}\bigr)\bigr)\bigr]\bigr\rrvert
\leq \E [K_{T}].%
\]
As a result,
\begin{eqnarray*}
\E\sup_{\psi\in\widetilde{\Psi}}\bigl\llvert \mathbb {G}_{n}\bigl(M(\psi)\bigr)1
\bigl\{ K_{T}>a(\sigma)\sqrt {n}\bigr\}\bigr\rrvert & \leq&\sqrt{n}\cdot
\E \bigl[ \bigl(K_{T}+\E[K_{T}]\bigr)1\bigl\{
K_{T}>a(\sigma )\sqrt{n}\bigr\} \bigr]
\\
& \leq& 2\sqrt{n}\cdot\E \bigl[ K_{T}\cdot1\bigl\{K_{T}>a(
\sigma)\sqrt {n}\bigr\} \bigr],
\end{eqnarray*}
where we set
\[
a(\sigma)=\frac{J(\sigma)}{\sqrt{\log[1+N(\sigma, \widetilde\Psi
, \rho)]}}.
\]
Under the condition $\| \mathfrak{M} \|_{\rho,T}\leq\Theta$ a.s. one
can prove that
%
\begin{equation}
\label{expineqMphMpsi} \P \Bigl(\mathop{\sup_{\psi,\phi\in\widetilde\Psi,}}_{ \rho
(\psi,\phi
) \leq\delta}
\sup_{t\in[0,T]}\bigl|M_{t}(\psi)-M_{t}(\phi)\bigr|>x \Bigr)\leq
2e^{-x^{2}/CJ^{2}(\delta)}
\end{equation}
for all $x>0,$ where
$ C $ is a universal constant depending only on $\Theta$.
Inequality~\eqref{expineqMphMpsi} implies
\[
\E \bigl[ K_{T}\cdot1\bigl\{K_{T}>a(\sigma)\sqrt{n}\bigr\}
\bigr]\leq 2a(\sigma )\sqrt{n} e^{-na^{2}(\sigma)/CJ^{2}(\sigma)} +2\int_{a(\sigma)\sqrt{n}}^{\infty}e^{-x^{2}/CJ^{2}(\sigma)}
\,dx.
\]
Fix an integer $q_{0}$ such that $\sigma\leq2^{-q_{0}}\leq2\sigma.$
For each natural number $ q>q_{0}, $ there exists a nested sequence
of partitions $\widetilde{\Psi}=\bigcup_{i=1}^{N_{q}}\widetilde{\Psi}_{qi}$ of $\widetilde\Psi$ into $ N_{q}
$ disjoint subsets such that
$\rho(\psi,\phi)\leq2^{-q}$ for any $\psi,\phi\in\widetilde
{\Psi
}_{qi}$ and $N_{q}\leq N(2^{-q+1},\widetilde\Psi,\rho).$ Denote
\[
\Delta_{qi}=\sup_{\psi,\phi\in\widetilde{\Psi}_{qi}}\sup_{t\in
[0,T]}\bigl\llvert
M_{t}(\psi)-M_{t}(\phi)\bigr\rrvert ,
\]
and then \eqref{expineqMphMpsi} implies
\[
\E\bigl[\Delta^{2}_{qi}\bigr]\leq2\int_{0}^{\infty}
xe^{-x^{2}/CJ^{2}(2^{-q})} \,dx=CJ^{2}\bigl(2^{-q}\bigr).
\]
Choose for each $q\geq q_{0}$ a fixed element $ \psi_{qi}$ from each
partioning set $\widetilde{\Psi}_{qi},$ and set
\[
\Pi_{q} \bigl[ Z\bigl(M(\psi)\bigr) \bigr] =Z(\psi_{qi}),\qquad
\Delta_{q} \bigl[ Z\bigl(M(\psi)\bigr) \bigr] =\Delta_{qi}
\qquad\mbox{if }\psi\in\widetilde{\Psi}_{qi}.
\]
Then $\Pi_{q} [ Z(M(\psi)) ] $ and $\Delta_{q} [
Z(M(\psi
)) ] $ run
through a set of $N_{q}$ functions if $\psi$ runs through $
\widetilde
{\Psi}. $
Define for each fixed $ n $ and $ q\geq q_{0} $ numbers and
indicator functions
\begin{eqnarray*}
&& a_{q}=J\bigl(2^{-q}\bigr)/\sqrt{\log[1+N_{q+1}]},
\\
&& A_{q-1}\bigl[Z\bigl(M(\psi)\bigr)\bigr]\\
&&\qquad=\mathbf{1}\bigl\{
\Delta_{q_{0}}\bigl[Z\bigl(M(\psi)\bigr)\bigr]\leq \sqrt
{n}a_{q_{0}},\ldots, \Delta_{q-1}\bigl[Z\bigl(M(\psi)\bigr)
\bigr]\leq\sqrt{n}a_{q-1}\bigr\},
\\
&& B_{q}\bigl[Z\bigl(M(\psi)\bigr)\bigr]\\
&&\qquad=\mathbf{1}\bigl\{
\Delta_{q_{0}}\bigl[Z\bigl(M(\psi)\bigr)\bigr]\leq \sqrt
{n}a_{q_{0}},\ldots, \Delta_{q-1}\bigl[Z\bigl(M(\psi)\bigr)
\bigr]\leq\sqrt{n}a_{q-1},
\\
&&\hspace*{202pt}\Delta_{q}\bigl[Z\bigl(M(\psi)\bigr)\bigr]>\sqrt{n}a_{q}
\bigr\}.
\end{eqnarray*}
Now decompose
\begin{eqnarray*}
\label{decompchain} &&Z\bigl(M(\psi)\bigr)-\Pi_{q_{0}}\bigl[Z\bigl(M(\psi)
\bigr)\bigr]\\
&&\qquad=\sum_{q=q_{0}+1}^{\infty} \bigl(Z\bigl(M(
\psi )\bigr)-\Pi_{q}\bigl[Z\bigl(M(\psi)\bigr)\bigr]\bigr)
B_{q}\bigl[Z\bigl(M(\psi)\bigr)\bigr]
\\
&&\qquad\quad{}+\sum_{q=q_{0}+1}^{\infty} \bigl(\Pi_{q}
\bigl[Z\bigl(M(\psi)\bigr)\bigr]-\Pi_{q-1}\bigl[Z\bigl(M(\psi)\bigr)
\bigr]\bigr) A_{q-1}\bigl[Z\bigl(M(\psi)\bigr)\bigr].
\end{eqnarray*}
We observe that either all of the $B_{q}[Z(M(\psi))]
$ are zero, in which case the $A_{q-1}[Z(M(\psi))]$ are $1$, or
alternatively, $ B_{q_{1}}[Z(M(\psi))]=1 $ for some $q_{1} > q_{0}$
(and zero for all other $q$), in which case $A_{q}[Z(M(\psi))]=1$
for $q < q_{1}$ and $A_{q}[Z(M(\psi))]=0$ for $q \geq q_{1}.$ Our
construction of partitions and choice of $ q_{0} $ also ensure that
\[
a(\sigma)=\frac{J(\sigma)}{\sqrt{\log[1+N(\sigma, \widetilde\Psi
, \rho
)]}}\leq\frac{J(2^{-q_{0}})}{\sqrt{\log[1+N(2^{-q_{0}-1},
\widetilde
\Psi, \rho)]}}\leq a_{q_{0}},
\]
whence $ A_{q_{0}}[Z(M(\psi))]=1. $
Next we apply the empirical process $ \mathbb{G}_{n} $\break to both series
on the right-hand side of separately, take absolute values, and\break
next
take suprema over $ \psi\in\widetilde{\Psi}. $ Because the
partitions are nested,\break $ \Delta_{q}[Z(M(\psi))] B_{q}[Z(M(\psi
))]\leq
\Delta_{q-1}[Z(M(\psi))] B_{q}[Z(M(\psi))]\leq\sqrt{n}a_{q-1}.$ The
last inequality holds if $B_{q}[Z(M(\psi))]=0$ and also if $
B_{q}[Z(M(\psi))]=1$ by definition. Furthermore, as $B_{q}[Z(M(\psi
))]$ is indicator of the event $\Delta_{q}[Z(M(\psi))]>\sqrt{n}
a_{q},$ it follows
\begin{eqnarray*}
\sqrt{n}a_{q}\cdot\E \bigl[\Delta_{q}\bigl[Z\bigl(M(
\psi)\bigr)\bigr] B_{q}\bigl[Z\bigl(M(\psi )\bigr)\bigr] \bigr] &\leq&\E
\bigl[\bigl(\Delta_{q}\bigl[Z\bigl(M(\psi)\bigr)\bigr]
\bigr)^{2} B_{q}\bigl[Z\bigl(M(\psi)\bigr)\bigr]\bigr]\\
&\leq&
J^{2}\bigl(2^{-q}\bigr)
\end{eqnarray*}
by the choice of $\Delta_{q}[Z(M(\psi))].$
Because $|\mathbb{G}_{n}[Z(M(\psi))]|\leq\mathbb{G}_{n}[Z']+2\sqrt {n}\cdot\E[Z']$ if $ |Z(M(\psi))|\leq Z',$ we obtain by the triangle
inequality and Lemma~\ref{expineqfinit} that the quantity
\[
\E\Biggl\llVert \sum_{q=q_{0}+1}^{\infty}
\mathbb{G}_{n}\bigl(Z\bigl(M(\psi)\bigr)-\Pi_{q}\bigl[Z
\bigl(M(\psi)\bigr)\bigr]\bigr)B_{q}\bigl[Z\bigl(M(\psi)\bigr)\bigr]
\Biggr\rrVert_{\widetilde{\Psi}}
\]
is bounded by
\begin{eqnarray*}
&&\sum_{q=q_{0}+1}\E\bigl\| \mathbb{G}_{n}
\Delta_{q}\bigl[Z\bigl(M(\psi )\bigr)\bigr]B_{q}\bigl[Z
\bigl(M(\psi )\bigr)\bigr] \bigr\|_{\widetilde{\Psi}}
\\
&&\quad{} +2\sqrt{n}\sum_{q=q_{0}+1}^{\infty}\bigl\| \E\bigl\{
\Delta_{q}\bigl[Z\bigl(M(\psi)\bigr)\bigr] B_{q}\bigl[Z
\bigl(M(\psi)\bigr)\bigr] \bigr\} \bigr\|_{\widetilde{\Psi}}
\\
&&\qquad\lesssim\sum_{q=q_{0}+1}^{\infty} \biggl[
a_{q-1}\log[1+N_{q}] +CJ\bigl(2^{-q}\bigr)\sqrt{
\log[1+ N_{q}]} +\frac{J^{2}(2^{-q})}{a_{q}} \biggr].
\end{eqnarray*}
In view of the definition of $ a_{q}, $ the series on the right can
be bounded by a multiple of the series $ \sum_{q=q_{0}+1}^{\infty}
J(2^{-q}) \sqrt{\log[1+N_{q}]}. $ To establish a similar bound for
the second part of equation \eqref{decompchain}, note that there are
at most $ N_{q} $ differences $ \Pi_{q}[Z(M(\psi))]-\Pi_{q-1}[Z(M(\psi))] $ and at most $ N_{q-1} $
indicator functions $
A_{q-1}[Z(M(\psi))]. $ Because the partitions are nested, $ (\Pi_{q}[Z(M(\psi))]-\Pi_{q-1}[Z(M(\psi))])A_{q-1}[Z(M(\psi))] $
is bounded by $ \Delta_{q-1}[Z(M(\psi))]A_{q-1}\times\break [Z(M(\psi))]\leq
\sqrt {n} a_{q-1}. $ Moreover, $ \E[\Pi_{q}[Z(M(\psi))]-\Pi_{q-1}[Z(M(\psi
))]]^{2}\leq\break CJ^{2}(2^{-q}). $ Hence
\begin{eqnarray*}
&&\Biggl\llVert \sum_{q_{0}+1}^{\infty}
\mathbb{G}_{n}\bigl(\Pi_{q}\bigl[Z\bigl(M(\psi )\bigr)
\bigr]-\Pi_{q-1}\bigl[Z\bigl(M(\psi)\bigr)\bigr]\bigr) A_{q-1}
\bigl[Z\bigl(M(\psi)\bigr)\bigr] \Biggr\rrVert_{\widetilde{\Psi}}
\\
&&\qquad\leq\sum_{q=q_{0}+1}^{\infty} \bigl[
a_{q-1}\log(1+N_{q}) +CJ\bigl(2^{-q}\bigr)\sqrt{
\log[1+N_{q}]} \bigr].
\end{eqnarray*}
Again this is bounded above by a multiple of the series $ \sum_{q=q_{0}+1}^{\infty}J(2^{-q})\times \sqrt{\log[1+N_{q}]}. $ To conclude the
proof it suffices to consider the terms\break $ \Pi_{q_{0}}[Z(M(\psi))]. $
Because $ |\Pi_{q_{0}}[Z(M(\psi))]|\leq K_{T} \leq a(\delta)\sqrt {n}\leq a_{q_{0}}\sqrt{n} $ and
\[
\E\bigl(\Pi_{q_{0}}\bigl[Z\bigl(M(\psi)\bigr)\bigr]\bigr)^{2}
\leq\E \Bigl[\sup_{t\in
[0,T]}\bigl(M_{t}(\psi_{q_{0}i})-M_{t}
\bigl(\psi^{*}\bigr)\bigr) \Bigr]^{2}\leq4\Theta^{2}
\sigma^{2}
\]
by the Burkholder--Davis--Gundy inequality, we have
\[
\E\bigl\| \mathbb{G}_{n}\Pi_{q_{0}}\bigl[Z\bigl(M(\psi)\bigr)\bigr]
\bigr\|_{\widetilde{\Psi
}}\lesssim a_{q_{0}}\log[1+N_{q_{0}}]+\sigma
\sqrt{\log[1+N_{q_{0}}]}.
\]
By the choice of $ q_{0}, $ this is bounded by a multiple of the
first few items of the series
\[
\sum_{q=q_{0}+1}^{\infty}J\bigl(2^{-q}
\bigr)\sqrt{\log[1+N_{q}]}.
\]
\upqed\end{pf}
%

\subsection{\texorpdfstring{Proof of Proposition \protect\ref{selfbound}}
{Proof of Proposition 6.1}}
The proof can be routinely carried out along with lines of \citet{BGR}.

\subsection{\texorpdfstring{Proof of Proposition \protect\ref{EVexpbound}}
{Proof of Proposition 6.3}}
\label{proofevexp}
In order to prove Proposition~\ref{EVexpbound} we need the following lemma.
%
\begin{lem}
\label{expineqfixedpsi}
Denote
\[
\mathcal{Q}\bigl(\psi,\psi'\bigr)=\rho\bigl(\psi,\psi'
\bigr)\sqrt{\log\log\bigl(\rho^{2}\bigl(\psi ,\psi'\bigr)
\vee e^{2}\bigr)}
\]
for any $\psi,\psi'\in\Psi.$
There is $\varepsilon>0$ such that for any $\psi,\psi'\in\Psi$
and $
\psi^{*}\in\Psi^{*}$, it holds
\[
\E \biggl\{\exp \biggl(\theta \biggl[\frac{\sqrt{n}\cdot(\Delta_{n}(\psi
)-\Delta_{n}(\psi'))}{\mathcal{Q}(\psi,\psi^{*})\cdot\mathcal
{Q}(\psi
,\psi')} \biggr] \biggr)-1 \biggr
\}\leq C\theta^{2}
\]
for some constant $C>0,$ provided $|\theta|\leq\varepsilon$.
\end{lem}
\begin{pf}
Without loss of generality, we may, and do, assume that \mbox{$\Theta=1$}.
Fix a martingale $M^{*}=M(\psi^{*})\in\mathcal{M}^{*}.$
Since $ Z(M^{*})=\E[Z(M^{*})] $\vadjust{\goodbreak} almost surely, we have for arbitrary
$M=M(\psi),M'=M(\psi')\in\mathcal{M}$
\begin{eqnarray*}
&&V_{n}(M)-V_{n}\bigl(M'\bigr)\\
&&\qquad=
\frac{1}{n(n-1)}\sum_{1\leq i<j\leq
n} \bigl(
\widetilde{Z}^{(i)}(M)-\widetilde{Z}^{(j)}(M) \bigr)^{2}
\\
&&\qquad\quad{} -\frac{1}{%
n(n-1)}\sum_{1\leq i<j\leq n} \bigl(
\widetilde{Z}^{(i)}\bigl(M'\bigr)-%
\widetilde{Z}^{(j)}\bigl(M^{\prime}\bigr) \bigr)^{2}
\\
&&\qquad=\frac{1}{n(n-1)}\sum_{1\leq i<j\leq n} \bigl(
\widetilde{Z}%
^{(i)}(M)-\widetilde{Z}^{(i)}
\bigl(M'\bigr)-\widetilde{Z}^{(j)}(M)+\widetilde{%
Z}^{(j)}\bigl(M'\bigr) \bigr)
\\
&&\qquad\quad{} \times \bigl( \widetilde{Z}^{(i)}(M)-\widetilde{Z}%
^{(j)}(M)+
\widetilde{Z}^{(i)}\bigl(M'\bigr)-\widetilde{Z}^{(j)}
\bigl(M'\bigr) \bigr)
\\
&&\qquad=\frac{1}{n(n-1)}\sum_{1\leq i<j\leq n} \bigl(
Z^{(i)}(M)-Z^{(i)}\bigl(M^{\prime}\bigr)-Z^{(j)}(M)+Z^{(j)}
\bigl(M^{\prime}\bigr) \bigr)
\\
&&\hspace*{76pt}\qquad\quad{}\times \bigl( \widetilde{Z}^{(i)}(M)-\widetilde{Z}^{(j)}(M)+
\widetilde{Z}%
^{(i)}\bigl(M^{\prime}\bigr)-
\widetilde{Z}^{(j)}\bigl(M^{\prime}\bigr) \bigr)
\end{eqnarray*}
with $\widetilde{Z}^{(i)}=Z^{(i)}(M)-Z^{(i)}(M^{*}),$ $i=1,\ldots
,n.$ Set
\[
\xi_{i} =Z^{(i)}(M)-Z^{(i)}\bigl(M^{\prime}
\bigr), \qquad\zeta_{i} =\widetilde{Z}^{(i)}(M)+\widetilde{Z}^{(i)}
\bigl(M^{\prime}\bigr),\qquad  i=1,\ldots, n,
\]
then
\begin{eqnarray*}
V_{n}(M)-V_{n}\bigl(M^{\prime}\bigr) &=&
\frac{1}{n(n-1)}\sum_{1\leq
i<j\leq
n} ( \xi_{i}-
\xi_{j} ) ( \zeta_{i}-\zeta_{j} )
\\
&=&\frac{1}{n(n-1)}\sum_{1\leq i<j\leq n}\xi_{i}
\zeta_{i}-\frac
{1}{%
n(n-1)}\sum_{1\leq i<j\leq n}
\xi_{i}\zeta_{j}
\\
&&{}-\frac{1}{n(n-1)}\sum_{1\leq i<j\leq n}\xi_{j}
\zeta_{i}+\frac
{1}{%
n(n-1)}\sum_{1\leq i<j\leq n}
\xi_{j}\zeta_{j}
\\
&=&\frac{2}{n}\sum_{i=1}^{n}
\xi_{i}\zeta_{i}-\frac
{1}{n(n-1)}\sum
_{i\neq j}\xi_{i}\zeta_{j}.
\end{eqnarray*}
Hence
%
\begin{eqnarray}\label{vnrepr}
&&
V_{n}(M)-V(M)-\bigl(V_{n}\bigl(M'
\bigr)-V\bigl(M'\bigr)\bigr)
\nonumber
\\[-8pt]
\\[-8pt]
\nonumber
&&\qquad=\frac{2}{n}\sum
_{i=1}^{n}\bigl(\xi_{i}
\zeta_{i}-\E[\xi_{i}\zeta_{i}]\bigr)
-\frac{1}{n(n-1)}\sum_{i\neq j}\bigl(
\xi_{i}\zeta_{j}-\E[\xi_{i}\zeta_{j}]
\bigr).
\end{eqnarray}
Note that $(\xi_{1},\zeta_{1}),\ldots, (\xi_{n},\zeta_{n})$ is a family
of i.i.d. random two-dimensional vectors such that
\[
|\xi_{i}|\leq\sup_{t\in[0,T]}\bigl|M^{(i)}_{t}-M'^{(i)}_{t}\bigr|,\qquad
i=1,\ldots, n
\]
and
\[
|\zeta_{i}|\leq2\sup_{t\in[0,T]}\bigl|M^{(i)}_{t}-M^{*(i)}_{t}\bigl|,\qquad
i=1,\ldots, n.
\]
Lemma~\ref{expineqlm} implies that for any $x>0$,
\[
\P \biggl( \frac{|\xi_{i}| }{\sqrt{%
\langle M^{(i)}-M'^{(i)} \rangle_{T}\log\log(
\langle
M^{(i)}-M'^{(i)} \rangle_{T}\vee e^{2})}}\geq x \biggr) \leq C(\alpha)e^{-\alpha x^{2}}
\]
and
\[
\P \biggl( \frac{|\zeta_{i}| }{\sqrt{%
\langle M^{(i)}-M^{*(i)} \rangle_{T}\log\log(
\langle
M^{(i)}-M^{*(i)} \rangle_{T}\vee e^{2})}}\geq x \biggr) \leq C(\alpha)e^{-\alpha x^{2}/4}.
\]
As a result,
\[
\P \biggl( \frac{|\xi_{i}| }{\rho(\psi,\psi')\sqrt{%
\log\log(\rho^{2}(\psi,\psi')\vee e^{2})}}\geq x \biggr) \leq C(\alpha )e^{-\alpha x^{2}}
\]
and
\[
\P \biggl( \frac{|\zeta_{i}| }{\rho(\psi,\psi^{*})\sqrt{%
\log\log(\rho^{2}(\psi,\psi^{*})\vee e^{2})}}\geq x \biggr) \leq C(\alpha)e^{-\alpha x^{2}/4}
\]
for $i=1,\ldots,n.$ Using representation \eqref{vnrepr}, we get
\begin{eqnarray*}
&&\frac{\sqrt{n}(\Delta_{n}(\psi)-\Delta_{n}(\psi'))}{\mathcal
{R}(\psi
,\psi^{*})\cdot\mathcal{R}(\psi,\psi')}\\
&&\qquad=\frac{2}{\sqrt{n}}\sum_{i=1}^{n}
\bigl(\widetilde\xi_{i}\widetilde\zeta_{i}-\E[\widetilde
\xi_{i}\widetilde\zeta_{i}]\bigr) -\frac{1}{\sqrt{n}(n-1)}\sum
_{i<j}\bigl(\widetilde\xi_{i}\widetilde
\zeta_{j}-\E[\widetilde\xi_{i}\widetilde
\zeta_{j}]\bigr)
\\
&&\qquad\quad{}-\frac{1}{\sqrt{n}(n-1)}\sum_{j<i}\bigl(\widetilde
\xi_{i}\widetilde \zeta_{j}-\E[\widetilde\xi_{i}
\widetilde\zeta_{j}]\bigr)\\
& =& T_{1,n}+T_{2,n}+T_{3,n},
\end{eqnarray*}
where the ``normalized'' random variables
\begin{eqnarray*}
\widetilde\xi_{i}&=&\frac{\xi_{i} }{\rho(\psi,\psi')\sqrt{%
\log\log(\rho^{2}(\psi,\psi')\vee e^{2})}},
\\
\widetilde\zeta_{i}&=&\frac{\zeta_{i} }{\rho(\psi,\psi^{*})\sqrt {%
\log\log(\rho^{2}(\psi,\psi^{*})\vee e^{2})}}
\end{eqnarray*}
satisfy
%
\begin{equation}
\label{expineqxizeta} \P \bigl( |\widetilde\zeta_{i}|\vee|\widetilde
\xi_{i}|\geq x \bigr) \leq C(\alpha)e^{-\alpha x^{2}/4},\qquad i=1,\ldots,n.
\end{equation}
The inequalities in \eqref{expineqxizeta} immediately imply
\begin{eqnarray*}
\P\bigl(|\widetilde\xi_{i}\widetilde\zeta_{i}|>x\bigr)&\leq&\P\bigl(|
\widetilde \xi_{i}|^{2}+|\widetilde\zeta_{i}|^{2}>2x
\bigr)\leq\P\bigl(|\widetilde\xi_{i}|>\sqrt{2x}\bigr)+\P\bigl(|\widetilde
\zeta_{i}|>\sqrt{2x}\bigr)\\
&\leq& 2C(\alpha )\exp (-\alpha x/2).
\end{eqnarray*}
Consider first the term $T_{1,n}.$ For any $\theta\in\mathbb{R}$
we have
%
\begin{equation}
\label{expT2} \E\bigl[\exp(\theta T_{1,n})\bigr]=\prod
_{i=1}^{n}\E\bigl[\exp\bigl(\theta\bigl(\widetilde
\xi_{i}\widetilde\zeta_{i}-\E[\widetilde\xi_{i}
\widetilde\zeta_{i}]\bigr)/\sqrt{n}\bigr)\bigr].
\end{equation}
Since the random variables $\widetilde\xi_{i}\widetilde\zeta_{i}-\E
[\widetilde\xi_{i}\widetilde\zeta_{i}],$
$i=1,\ldots,n,$ possess finite moments of any order and have zero
mean, it holds
\[
\log\E \bigl[\exp\bigl(\varepsilon\bigl(\widetilde\xi_{i}\widetilde
\zeta_{i}-\E [\widetilde\xi_{i}\widetilde
\zeta_{i}]\bigr)\bigr) \bigr]= \tfrac{1}{2}\sigma^{2}
\varepsilon^{2}+o\bigl(\varepsilon^{2}\bigr),\qquad  i=1,\ldots,n
\]
as $\varepsilon\to0,$ where $\sigma^{2}=\E (\widetilde\xi_{i}\widetilde\zeta_{i}-\E[\widetilde\xi_{i}\widetilde\zeta_{i}]
)^{2}.$ Hence the inequality
%
\begin{equation}
\label{momgen} \E \bigl[\exp\bigl(\varepsilon\bigl(\widetilde\xi_{i}
\widetilde\zeta_{i}-\E [\widetilde\xi_{i}\widetilde
\zeta_{i}]\bigr)\bigr) \bigr]\leq e^{C_{1}\varepsilon^{2}}
\end{equation}
holds for sufficiently small $\varepsilon$ and any $C_{1}>\sigma^{2}/2.$ Combining \eqref{expT2} with
\eqref{momgen}, we get for all $n\in\mathbb{N}$ and sufficiently
small $\theta>0$,
\[
\E\bigl[\exp(\theta T_{1,n})-1\bigr]\leq e^{C_{1}\theta^{2}}-1\leq
C_{2}\theta^{2}.
\]
Turn now to the terms $T_{2,n}$ and $T_{3,n}.$ We need the
following proposition to estimate~$T_{2,n}$
and $T_{3,n}.$
%
\begin{prop}
\label{lemmomineq}
Let $(X_{1},Y_{1}),\ldots,(X_{n},Y_{n})$ be a sequence of i.i.d.
centered random vectors in $\mathbb{R}^{2}$
such that $\E|X_{i}|^{p}<\infty$ and $\E|Y_{i}|^{p}<\infty$ for
all $i=1,\ldots,n, $ and some $p\geq2.$ Then
%
\begin{eqnarray}\label{momineq}
&&\E\biggl\llvert \sum_{1\leq i<j\leq n}X_{i}Y_{j}
\biggr\rrvert^{p}\nonumber\\
&&\qquad\leq C^{p}\max \Biggl\{\sum
_{1\leq i<j\leq n}\E|X_{i}|^{p}\E|Y_{j}|^{p},
\sum_{i=1}^{n-1}\E |X_{i}|^{p}
\Biggl(\sum_{j=i+1}^{n}\E|Y_{j}|^{2}
\Biggr)^{p/2},
\\
 &&\hspace*{52pt} \qquad \sum_{j=2}^{n}
\E|Y_{j}|^{p} \Biggl(\sum_{i=1}^{j-1}
\E |X_{i}|^{2} \Biggr)^{p/2}, \biggl(\sum
_{1\leq i<j\leq n}\E |X_{i}|^{2}\E
|Y_{j}|^{2} \biggr)^{p/2} \Biggr\}\nonumber
\end{eqnarray}
for some constant $C>0$ not depending on $p$.
\end{prop}
\begin{pf}
Denote $Q_{n}=\sum_{1\leq i<j\leq n}X_{i}Y_{j}$ and\vspace*{1pt}
\[
V_{j}=\sum_{i=1}^{j-1}X_{i}Y_{j},\qquad
j=2,\ldots,n.\vspace*{-1pt}
\]
It is clear that $T_{2,n}=\sum_{j=2}^{n}V_{j}$ and $(V_{j},
j=2,\ldots,n)$ is a forward martingale-difference sequence (see the
\hyperref[app]{Appendix} for definition) with respect to $\sigma$-algebras $
\mathcal
{F}_{j}=\sigma((X_{1},Y_{1}),\ldots, (X_{j},Y_{j})), j=2,\ldots,
n.$ By the martingale Rosenthal inequality (see Proposition~\ref
{burkhold} in the \hyperref[app]{Appendix}),\vspace*{-1pt}
\[
\E\bigl[|Q_{n}|^{p}\bigr]\leq B(p/\log p)\max \Biggl\{\sum
_{j=2}^{n}\E|V_{j}|^{p},
\E \Biggl(\sum_{j=2}^{n}\E \bigl[
V_{j}^{2}\rrvert \mathcal {F}_{j-1} \bigr]
\Biggr)^{p/2} \Biggr\}\vspace*{-1pt}
\]
and\vspace*{-1pt}
%
\begin{equation}
\label{t2n1} \E|V_{j}|^{p}\leq B(p/\log p) \cdot
\E|Y_{j}|^{p}\max \Biggl\{\sum
_{i=1}^{j-1}\E|X_{i}|^{p},
\Biggl(\sum_{i=1}^{j-1}\E
|X_{i}|^{2} \Biggr)^{p/2} \Biggr\}\vspace*{-1pt}
\end{equation}
for all $j=2,\ldots,n.$
Then\vspace*{-1pt}
\begin{eqnarray*}
&&\E \Biggl(\sum_{j=2}^{n}\E \bigl[
V^{2}_{j}\rrvert \mathcal {F}_{j-1} \bigr]
\Biggr)^{p/2}
\\[-1pt]
&&\qquad=\E \Biggl(\sum_{1\leq i<j\leq n}|X_{i}|^{2}
\E|Y_{j}|^{2}+2\sum_{j=3}^{n}
\sum_{1\leq k<l\leq j-1}X_{k}X_{l}
\E|Y_{j}|^{2} \Biggr)^{p/2}
\\[-1pt]
&&\qquad\leq2^{p/2-1} \E \biggl(\sum_{1\leq i<j\leq n}|X_{i}|^{2}
\E |Y_{j}|^{2} \biggr)^{p/2}\\
&&\qquad\quad{} +2^{p-1}\E
\Biggl\llvert \sum_{1\leq k<l\leq
n-1}X_{k}X_{l}
\sum_{j=l+1}^{n}\E|Y_{j}|^{2}
\Biggr\rrvert^{p/2}.\vspace*{-1pt}
\end{eqnarray*}
By the Rosenthal inequality,\vspace*{-1pt}
%
\begin{eqnarray}
\label{t2n2} &&\E \biggl(\sum_{1\leq i<j\leq n}|X_{i}|^{2}
\E|Y_{j}|^{2} \biggr)^{p/2}\nonumber\\
&&\qquad=\E \Biggl(\sum
_{i=1}^{n-1}|X_{i}|^{2}\sum
_{j=i+1}^{n}\E |Y_{j}|^{2}
\Biggr)^{p/2}
\nonumber
\\[-10pt]
\\[-10pt]
\nonumber
&&\qquad\leq B(p/2)\log^{-1}(p/2)\max \Biggl\{\sum_{i=1}^{n-1}
\E|X_{i}|^{p} \Biggl[\sum_{j=i+1}^{n}
\E|Y_{j}|^{2} \Biggr]^{p/2},\\
&&\hspace*{154pt} \biggl(\sum
_{1\leq i<j\leq
n}\E |X_{i}|^{2}\E|Y_{j}|^{2}
\biggr)^{p/2} \Biggr\}.\nonumber
\end{eqnarray}
Using the Jensen inequality, we get for $2\leq p<4$
\begin{eqnarray*}
&&\E\Biggl\llvert \sum_{1\leq k<l\leq n-1}X_{k}X_{l}
\sum_{j=l+1}^{n}\E |Y_{j}|^{2}
\Biggr\rrvert^{p/2}\\[-2pt]
&&\qquad\leq \Biggl(\sum_{1\leq k<l\leq n-1}\E
|X_{k}X_{l}|^{2} \Biggl(\sum
_{j=l+1}^{n}\E|Y_{j}|^{2}
\Biggr)^{2} \Biggr)^{p/4}.
\end{eqnarray*}
Moreover,
%
\begin{eqnarray}
\label{t2n3} &&\Biggl(\sum_{1\leq k<l\leq n-1}\E|X_{k}X_{l}|^{2}
\Biggl(\sum_{j=l+1}^{n}\E
|Y_{j}|^{2} \Biggr)^{2} \Biggr)^{p/4}
\nonumber
\\[-9pt]
\\[-9pt]
\nonumber
&&\qquad\leq \biggl(\sum_{1\leq i<j\leq n}\E|X_{i}|^{2}
\E|Y_{j}|^{2} \biggr)^{p/2}.
\end{eqnarray}
Combining \eqref{t2n1}, \eqref{t2n2} and \eqref{t2n3}, we arrive at the
inequality \eqref{momineq}. Thus Lemma~\ref{lemmomineq} is proved for all $2\leq p <4.$ Suppose now that
the inequality \eqref{momineq} holds for $p\leq m-1$ with some
$m>4.$ Let us prove it for $p=m.$ It follows from the previous
steps, that we only need to obtain an upper bound for the term
\[
\E\Biggl\llvert \sum_{1\leq k<l\leq n-1}X_{k}X_{l}
\sum_{j=l+1}^{n}\E |Y_{j}|^{2}
\Biggr\rrvert^{m/2}.
\]
Our induction hypothesis gives that the quantity
\[
\E\Biggl\llvert \sum_{1\leq k<l\leq n-1}X_{k}X_{l}
\sum_{j=l+1}^{n}\E |Y_{j}|^{2}
\Biggr\rrvert^{m/2}
\]
is bounded by
%
\begin{eqnarray}
\label{t2n4}&& C^{m/2}\max \Biggl\{\sum_{1\leq k<l\leq n-1}
\E\Biggl\llvert X_{k}X_{l}\sum
_{j=l+1}^{n}\E|Y_{j}|%
^{2}
\Biggr\rrvert^{m/2},
\nonumber\\[-2pt]
&&\hspace*{33pt}\qquad\sum_{k=1}^{n-2}
\E|X_{k}|^{m/2} \Biggl(\sum_{l=k+1}^{n-1}|X_{l}|^{2}
\sum_{j=l+1}^{n}\E|Y_{j}|^{2}
\Biggr)^{m/4},
\nonumber
\\[-9pt]
\\[-9pt]
\nonumber
&&\hspace*{33pt}\qquad\sum_{l=2}^{n-1}|X_{l}|^{m/2}
\Biggl(\sum_{j=l+1}^{n}\E
|Y_{j}|^{2} \Biggr)^{m/2} \Biggl(\sum
_{k=1}^{l-1}\E|X_{k}|^{2}
\Biggr)^{m/4},
\\[-2pt]
&&\hspace*{37pt}\qquad\Biggl(\sum_{1\leq k<l\leq n-1}\E
\bigl[|X_{k}|^{2}|X_{l}|^{2} \bigr]
\Biggl(\sum_{j=l+1}^{n}\E |Y_{j}|^{2}
\Biggr)^{2} \Biggr)^{m/4} \Biggr\}.\nonumber
\end{eqnarray}
Let us consider, for example, the first term in the above maximum.
Using the inequality
%
\begin{equation}
\label{ineqU} \Biggl( \E\sum_{k=1}^{n}
\llvert U_{k}\rrvert^{p} \Biggr)^{2}\leq \max
\Biggl\{ \sum_{k=1}^{n}\E\llvert
U_{k}\rrvert^{2p}, \Biggl( \sum
_{k=1}^{n}\E\llvert U_{k}\rrvert
\Biggr)^{2p} \Biggr\}
\end{equation}
that holds for any $p>1$ and any sequence of independent r.v. $
U_{1},\ldots,U_{n}$ with $\E|U_{k}|^{2p}<\infty,$ we get
\begin{eqnarray*}
\sum_{1\leq k<l\leq n-1}\E\Biggl\llvert X_{k}X_{l}
\sum_{j=l+1}^{n}\E|Y_{j}|
^{2}\Biggr\rrvert^{m/2} & =&\sum
_{1\leq k<l\leq n-1}\E\llvert X_{k}%
X_{l}
\rrvert^{m/2} \Biggl( \sum_{j=l+1}^{n}
\E|Y_{j}|^{2} \Biggr)^{m/2}
\\
& \leq& \Biggl[ \sum_{i=1}^{n-1}\E\llvert
X_{i}\rrvert^{m/2} \Biggl( \sum_{j=i+1}^{n}
\E|Y_{j}|^{2} \Biggr)^{m/4} \Biggr]^{2}
\end{eqnarray*}
and
\begin{eqnarray*}
&&\Biggl[ \sum_{i=1}^{n-1}\E\llvert
X_{i}\rrvert^{m/2} \Biggl( \sum_{j=i+1}^{n}
\E|Y_{j}|^{2} \Biggr)^{m/4} \Biggr]^{2}\\
&&\qquad
\leq\max \Biggl\{ \sum_{i=1}^{n-1}\E\llvert
X_{i}\rrvert^{m} \Biggl( \sum_{j=i+1}^{n}
\E|Y_{j}|^{2} \Biggr)^{m/2},
 \biggl( \sum_{1\leq i<j\leq n}%
\E|X_{i}|^{2}\E|Y_{j}|^{2}
\biggr)^{m/2} \Biggr\}.
\end{eqnarray*}
To see that inequality \eqref{ineqU} holds, just note that the function
\[
h(t)=\log \Biggl[ \sum_{k=1}^{n}
\E|U_{k}|^{t} \Biggr]
\]
is convex in the domain $t>1.$ Due to convexity of $h(t),$ we have %
\[
\Biggl( \sum_{k=1}^{n}
\E|U_{k}|^{p} \Biggr)^{2p-1}\leq \biggl( \sum
_{k=1}%
^{n}\E|U_{k}|^{2p}
\biggr)^{p-1} \Biggl( \sum_{k=1}^{n}
\E |U_{k}| \Biggr)^{p}%
\]
for any $p>1.$ Hence%
\begin{eqnarray*}
\Biggl( \sum_{k=1}^{n}\E
U_{k}^{p} \Biggr)^{2} & \leq& \biggl( \sum
_{k=1}%
^{n}\E|U_{k}|^{2p}
\biggr)^{{2(p-1)}/{(2p-1)}} \Biggl( \sum_{k=1}^{n}%
\E|U_{k}| \Biggr)^{{2p}/{(2p-1)}}
\\
& \leq&\max \Biggl\{ \sum_{k=1}^{n}
\E|U_{k}|^{2p}, \Biggl( \sum_{k=1}^{n}%
\E|U_{k}| \Biggr)^{2p} \Biggr\} .
\end{eqnarray*}
Other terms on the right-hand side of \eqref{t2n4} can be handled in a
similar way.
\end{pf}
Let us proceed with estimating the term $T_{2,n}.$ Without loss of
generality we may assume that $\E[\widetilde\xi]=\E[\widetilde
\zeta
]=0.$ Note that
for any natural $p>0$,
\begin{eqnarray*}
\E\bigl[|\widetilde\xi|^{p}\bigr] &\leq&2p C(\alpha)\int
_{0}^{\infty
}x^{p-1}\exp\bigl(-\alpha
x^{2}\bigr)\,dx\\
&=&\frac{2p C(\alpha)}{(2\alpha)^{p/2}}\int_{0}^{\infty
}y^{p-1}
\exp\bigl(-y^{2}/2\bigr)\,dy
\\
&\leq&\frac{p\sqrt{2\pi}C(\alpha)}{(2\alpha)^{p/2}}\E\bigl[|Z|^{p}\bigr],
\end{eqnarray*}
where $Z\sim N(0,1).$ Similarly
\[
\E\bigl[|\widetilde\zeta|^{p}\bigr]\leq\frac{2^{p/2}p\sqrt{2\pi}C(\alpha
)}{\alpha^{p/2}}\E
\bigl[|Z|^{p}\bigr].
\]
As a result, we get from Proposition~\ref{lemmomineq}
\begin{eqnarray*}
&&\E\bigl[|T_{2,n}|^{p}\bigr]\leq C^{p}\max \bigl
\{n^{1-p/2}(n-1)^{1-p}\E \bigl[|Z|^{2p}
\bigr],\\
&&\hspace*{98pt}n^{-p/2}(n-1)^{1-p/2}\E\bigl[|Z|^{p}
\bigr],(n-1)^{-p/2} \bigr\}
\end{eqnarray*}
for some constant $C>0$ and any $p>1.$ Hence for any $\theta\in
\mathbb{R}$,
%
\begin{eqnarray}
\label{expT2n} \E\bigl[\exp(\theta T_{2,n})-1\bigr]&=&\sum
_{k=2}^{\infty}\frac{\theta^{k}}{k!}\E \bigl[T_{2,n}^{k}
\bigr]
\nonumber\\
&\leq& \sum_{k=2}^{\infty}
\frac{|\theta|^{k}}{k!}\frac
{B_{1}^{k}}{(n-1)^{k/2}}\E\bigl[Z^{2k}\bigr]
\nonumber\\
&=&\E\bigl[\exp\bigl(B_{1}|\theta| Z^{2}/
\sqrt{n-1}\bigr)\bigr]-1-B_{1}|\theta|\E \bigl[Z^{2}\bigr]/
\sqrt{n-1}
\\
&=&\frac{1}{\sqrt{1-2B_{1}|\theta|/\sqrt{n-1}}}-1-B_{1}|\theta |/\sqrt{n-1}
\nonumber\\
&\leq&B_{2}\theta^{2},\nonumber
\end{eqnarray}
provided $B_{1}|\theta| /\sqrt{n-1}<1/2,$ where $B_{1}$ and $
B_{2}$ are two constants not depending on $k$ and $n.$
Analogously to \eqref{expT2n}, one can prove that
\[
\E\bigl[\exp(\theta T_{3,n})-1\bigr]\leq B_{3}
\theta^{2}
\]
for sufficiently small $|\theta|.$
Hence by the Cauchy--Schwarz inequality,
\begin{eqnarray*}
\E \bigl[e^{\theta(T_{1,n}+T_{2,n}+T_{3,n})}-1 \bigr]&\leq& \bigl[\E e^{2\theta T_{1,n}}
\bigr]^{1/2} \bigl[\E e^{4\theta T_{2,n}} \bigr]^{1/4} \bigl[\E
e^{4\theta T_{3,n}} \bigr]^{1/4}-1
\\
&\leq& B_{4}\theta^{2}
\end{eqnarray*}
for some constant $B_{4}>0.$
Lemma~\ref{expineqfixedpsi} is proved.
\end{pf}
Let us proceed with the proof of Proposition~\ref{EVexpbound}.
Let $ \{\widetilde\Psi^{m} \}_{m\in\mathbb{N}}$ be a
sequence of finite subsets of $\widetilde\Psi$ such that\vadjust{\goodbreak} $
\widetilde
\Psi^{m}\uparrow\widetilde\Psi$ as $m\to\infty.$ Introduce the
disjoint sets
\[
H_{p}=\bigl\{\psi\in\widetilde\Psi\dvtx  2^{-p-1}<\rho\bigl(\psi,
\psi^{*}\bigr) \leq 2^{-p}\bigr\}
\]
for any $p\in\mathbb{Z}.$ Without loss of generality we may assume
that $H_{p}$ are empty for $p<0.$
For every $m\in\mathbb{N},$ denote by $q(m,p)$ the smallest
integer such that $q(m,p)>p$ and that each of the closed balls with
centers in $\widetilde\Psi^{m}\cap H_{p}$ and $\rho$-radius $
2\cdot2^{-q(m,p)}$ contains exactly one point in $\widetilde\Psi^{m}\cap H_{p}.$ Then it is clear that
$\operatorname{Card}(\widetilde\Psi^{m}\cap H_{p})\leq
N(2^{-q(m,p)},\widetilde\Psi\cap H_{p},\rho).$ Next let us introduce
some mappings $\pi_{r}^{m,p}\dvtx \widetilde\Psi^{m}\cap H_{p}\to
\widetilde\Psi_{r}^{m,p},$ $p\leq r \leq q(m,p),$ defined by
\[
\pi_{r}^{m,p}=\lambda^{m,p}_{r}\circ
\lambda^{m,p}_{r+1}\circ\cdots \circ\lambda^{m,p}_{q(m,p)},
\]
where the sets $\widetilde\Psi_{r}^{m,p}\subset\widetilde\Psi^{m}\cap
H_{p}$ and the mappings $\lambda^{m,p}_{r}\dvtx \widetilde\Psi^{m}\cap
H_{p}\to\widetilde\Psi_{r}^{m,p}$
are specified in the following way. For $p\leq r<q(m,p),$ choose $
\widetilde\Psi_{r}^{m,p}$ and define $\lambda^{m,p}_{r}$ such that
they satisfy the following two conditions: $\operatorname
{Card}(\widetilde\Psi^{m,p}_{r})\leq N(2^{-r},\widetilde\Psi\cap
H_{p},\rho)$ and $\rho(\psi,\lambda^{m,p}_{r}(\psi))\leq2\cdot
2^{-r}$ for every $\psi\in\widetilde\Psi^{m}\cap H_{p}.$ For $
r=q(m,p),$ put $\widetilde\Psi^{m,p}_{q(m,p)}=\widetilde\Psi^{m}\cap
H_{p}$ and denote by $\lambda^{m,p}_{q(m,p)}$ the identical mapping
on $\widetilde\Psi^{m}\cap H_{p}.$ In terms of the mappings $\pi_{r}^{m,p}$ which have been introduced, we consider the chaining given
as follows: for every $n\in\mathbb{N}$ and $\psi\in\widetilde
\Psi
\cap H_{p}$,
\[
\bigl|\Delta_{n}(\psi)\bigr|\leq\sum_{r=p+1}^{q(m,p)}
\bigl|\Delta_{n}\bigl(\pi_{r}^{m,p}(\psi)\bigr)-
\Delta_{n}\bigl(\pi_{r-1}^{m,p}(\psi)\bigr)\bigr|+\bigl|
\Delta_{n}\bigl(\pi_{p}^{m,p}(\psi)\bigr)\bigr|.
\]
Since $\rho(\pi_{r}^{m,p}(\psi),\pi_{r-1}^{m,p}(\psi))/\rho(\psi
,\psi^{*})\leq2^{-r+p+1}$ and $\rho(\pi_{r}^{m,p}(\psi),\psi^{*})/ \rho
(\psi,\psi^{*})\leq2$
on $\widetilde\Psi^{m}\cap H_{p},$ it follows from Lemma~\ref
{expineqfixedpsi} and Lemma 8.2 in \citet{Kos} that
\begin{eqnarray*}
&&\E \biggl[\exp \biggl(\theta\sup_{\psi\in\widetilde\Psi^{m}\cap
H_{p}} \biggl\{\frac{\sqrt{n}|\Delta_{n}(\pi_{r}^{m,p}(\psi))-\Delta_{n}(\pi_{r-1}^{m,p}(\psi))|}{\mathcal{R}^{2}(\psi,\psi^{*})}
\biggr\} \biggr)-1 \biggr]
\\
&&\qquad\leq\E \biggl[\exp \biggl(\theta\sup_{\psi\in\widetilde\Psi
^{m}\cap
H_{p}} \biggl\{\frac{\mathcal{Q}(\pi_{r-1}^{m,p}(\psi),\pi_{r}^{m,p}(\psi
))\mathcal{Q}(\pi_{r}^{m,p}(\psi),\psi^{*})}{\mathcal{R}^{2}(\psi
,\psi^{*})}
\\
&&\hspace*{88pt}\qquad\quad{}\times\frac{\sqrt{n}|\Delta_{n}(\pi_{r}^{m,p}(\psi
))-\Delta_{n}(\pi_{r-1}^{m,p}(\psi))|}{\mathcal{Q}(\pi_{r-1}^{m,p}(\psi
),\pi_{r}^{m,p}(\psi))\mathcal{Q}(\pi_{r}^{m,p}(\psi),\psi^{*})} \biggr\} \biggr)-1 \biggr]
\\
&&\qquad\leq K p^{-2}4^{-r+p+1}\log\bigl(1+N\bigl(2^{-r},
\widetilde\Psi\cap H_{p},\rho\bigr)\bigr)
\end{eqnarray*}
for all $|\theta|\leq\varepsilon$, some $\delta>0$ and some
constant $K>0.$
Moreover note that $N(2^{-r},\widetilde\Psi\cap H_{p},\rho)\leq
N(2^{-r+p+1},\widetilde\Psi,\rho).$
Next
\begin{eqnarray*}
&&\E \biggl[\exp \biggl(\theta\sup_{\psi\in\widetilde\Psi^{m}\cap
H_{p}} \biggl\{\frac{|\sqrt{n}\cdot\Delta_{n}(\pi_{p}^{m,p}(\psi))|}{\mathcal
{R}^{2}(\psi,\psi^{*})}
\biggr\} \biggr)-1 \biggr]\\
&&\qquad\lesssim p^{-2}\log \bigl(1+N
\bigl(2^{1+p},\widetilde\Psi,\rho\bigr)\bigr).
\end{eqnarray*}
Finally, we get for any $P>0$,
\begin{eqnarray*}
& &\E \biggl[\exp \biggl(\theta\sup_{\psi\in\widetilde\Psi
^{m}\cap
(H_{1}\cup\cdots\cup H_{P})}\frac{|\sqrt{n}\cdot\Delta_{n}(\psi
)|}{\mathcal{R}^{2}(\psi,\psi^{*})} \biggr)-1
\biggr]
\\
&&\qquad\lesssim \sum_{p=1}^{P}p^{-2}
\sum_{r=p+1}^{q(m,p)} 4^{-r+p+1}\log
\bigl(1+N\bigl( 2^{-r+p+1}, \widetilde\Psi,\rho\bigr)\bigr)
\\
&&\qquad\lesssim\sum_{p=1}^{P}p^{-2}
\int_{0}^{1}\log\bigl(1+N(\sqrt {\varepsilon },
\widetilde\Psi,\rho)\bigr) \,d\varepsilon
\\
&&\qquad \lesssim \int_{0}^{1}\sqrt{\log\bigl(1+N(
\varepsilon,\widetilde\Psi ,\rho )\bigr)} \,d\varepsilon.
\end{eqnarray*}
The proof of Proposition~\ref{EVexpbound} is accomplished by letting
$m\to\infty$ and $P\to\infty.$

\begin{appendix}\label{app}

\section*{Appendix}
The following lemma is a straightforward generalization of Lemma 19.33
in \citet{VV}.
%
\begin{lemm}
\label{expineqfinit}
Let $ \mathcal{X} $ be a finite collection of bounded real valued
random variables defined on a common probability space $ (\Omega,
\mathcal{F}, \P), $ then
\[
\E\bigl\| \mathbb{G}_{n}[X] \bigr\|_{\mathcal{X}}\lesssim\frac{\sup_{X\in
\mathcal{X}}|X|}{\sqrt{n}}
\log\bigl(1+|\mathcal{X}|\bigr)+\max_{X\in
\mathcal
{X}}\sqrt{\E \bigl[|X|^{2} \bigr]}
\sqrt{\log\bigl(1+|\mathcal{X}|\bigr)},
\]
where $ \mathbb{G}_{n}[X]=\frac{1}{n}\sum_{j=1}^{n} (X^{(j)}-\E
[X]) $
and $ X^{(1)},\ldots, X^{(n)} $
are i.i.d. copies of~$X$.
\end{lemm}

Given a sequence of $\sigma$-algebras $(\mathcal{F}_{n}), n\geq
1$ on some probability space $(\Omega,\mathcal{F},\P)$,
we call a sequence of integrable r.v. $Y_{n}$ to be a forward
martingale-difference sequence w.r.t. $(\mathcal{F}_{n})$ if:
\begin{itemize}
\item$\mathcal{F}_{1}\subseteq\mathcal{F}_{2}\subseteq\cdots$;
\item$Y_{n}$ is $\mathcal{F}_{n}$-measurable;
\item$\E[Y_{n}|\mathcal{F}_{n-1}]=0$ a.s. for any $n\geq1$.
\end{itemize}
The following proposition can be found in \citet{Hit}.
%
\begin{propp}
\label{burkhold}
Let $(X_{k})$ be a forward martingale-difference sequence relative to
$\mathcal{F}_{k}$ such that
$\E|X_{k}|^{p}<\infty$ for some $p\geq2$ and $k=1,\ldots, n;$ then
\[
\E\Biggl\llvert \sum_{k=1}^{n}X_{k}
\Biggr\rrvert^{p}\leq B \bigl(k\log^{-1} k\bigr)\max \Biggl\{
\sum_{k=1}^{n}\E|X_{k}|^{p},
\E \Biggl[\sum_{k=1}^{n}\E
\bigl[X^{2}_{k}|\mathcal {F}_{k-1}\bigr]
\Biggr]^{p/2} \Biggr\}
\]
for some constant $B$ not depending on $k.$
\end{propp}
The next inequality can be found in \citet{DKL}.
%
\begin{lemm}
\label{expineqlm}
For any continuous local martingale $(M_{t})_{t\in\lbrack0,T]}$ with\break \mbox{$M_{0}=0$}%
\[
\P \biggl( \frac{\sup_{0\leq t\leq T}\llvert  M_{t}\rrvert
}{\sqrt{%
\langle M \rangle_{T}\log\log( \langle M
\rangle_{T}\vee e^{2})}}\geq x \biggr) \leq C(\alpha)e^{-\alpha x^{2}},
\]
where $\alpha$ is a real number in $(0,1/2)$ and $C(\alpha)$ is a positive
constant.
\end{lemm}
\end{appendix}

\section*{Acknowledgments}
I would like to thank John Schoenmakers,
Vladimir Spokoiny and Mikhail Urusov for remarks and helpful discussions.

%

%


\printaddresses

\end{document}